%% file: convex-chains.tex
\documentclass[11pt]{amsart}

\usepackage[a4paper]{geometry}
\usepackage[utf8x]{inputenc}
\usepackage[english]{babel}
\usepackage{graphicx}
\usepackage{lmodern}

\usepackage{enumerate}
\usepackage{xspace}

\usepackage[unicode]{hyperref}
\usepackage[abbrev,msc-links]{amsrefs}

\newcommand{\XX}{\mathbb{X}}
\newcommand{\ZZ}{\mathbb{Z}}
\newcommand{\NN}{\mathbb{N}}
\newcommand{\RR}{\mathbb{R}}

\newcommand{\CC}{\mathbb{C}}
\newcommand{\EE}{\mathbb{E}}
\newcommand{\PP}{\mathbb{P}}
\renewcommand{\L}{\XX}

\newcommand{\Ecal}{\mathcal{E}}
\newcommand{\bigO}{O}

\newcommand{\Sinai}{Sina{\u{\i{}}}\xspace}
\newcommand{\Moebius}{M\"{o}bius\xspace}
\newcommand{\Barany}{B\'ar\'any\xspace}
\newcommand{\Jarnik}{Jarn\'{\i}k\xspace}
\newcommand{\Zunic}{{\v{Z}}uni{\'c}\xspace}
\newcommand{\Erdos}{Erd{\"o}s\xspace}

\newcommand{\Z}{\ZZ}
\newcommand{\N}{\NN}

\newcommand{\Lstar}{{\L}}
\newcommand{\ds}{\displaystyle}
\newcommand{\PPbl}{\PP_{\beta,\lambda}}

\newcommand{\EEbl}{\EE_{\beta,\lambda}}

\newcommand{\Gammabl}{\Gamma_{\beta,\lambda}}
\newcommand{\sigmabl}{\sigma_{\beta,\lambda}}
\newcommand{\Lbl}{L_{\beta,\lambda}}
\newcommand{\Zbl}{Z(\beta,\lambda)}
\newcommand{\phibl}{\phi_{\beta,\lambda}}

\def\e{{\bf e}}
\def\c{{\bf c}}

\def\LL{\rm{Li}_2}
\def\LLL{\rm{Li}_3}

\renewcommand{\epsilon}{\varepsilon}

\newcommand{\indicator}[1]{\mathbf{1}_{\{#1\}}}

\newtheorem{lem}{Lemma}
\newtheorem{theorem}{Theorem}
\newtheorem*{theorem*}{Theorem}
\newtheorem{remark}{Remark}[theorem]
\newtheorem*{remark*}{Remark}
\newtheorem*{conjecture}{Conjecture}

\DeclareMathOperator{\Li}{Li}

\title{Lattice convex chains in the plane}

\author[J. Bureaux]{Julien Bureaux}
\address[J. Bureaux]{MODAL'X, Université Paris Ouest Nanterre La Défense, 200 avenue de la République, 92001 Nanterre}
\email{julien.bureaux@u-paris10.fr}

\author[N. Enriquez]{Nathana\"el Enriquez}
\address[N. Enriquez]{MODAL'X, Université Paris Ouest Nanterre La Défense, 200 avenue de la République, 92001 Nanterre}
\email{nathanel.enriquez@u-paris10.fr}
\address[N. Enriquez]{LPMA, Université Pierre et Marie Curie, 4 place Jussieu, 75005 Paris}

\begin{document}

\begin{abstract}
A detailed combinatorial analysis of planar lattice convex polygonal lines
is presented. This makes it possible to answer an open question of Vershik regarding the existence of a limit shape when the number of vertices is constrained. The method which is used emphasizes the connection of the combinatorial analysis with the zeros of the zeta function. It is shown how the Riemann Hypothesis leads to an asymptotic equivalent of the number of convex chains.
\end{abstract}

\maketitle

\section{Introduction}

In 1926, Jarn\'{\i}k found an equivalent of the maximal number of integral points
that a portion of length $n$ of the graph of a strictly convex function can interpolate. He obtained  an explicit constant times $n^{2/3}$. This work was at the origin of many works of Diophantine analysis, and we
refer the reader to the papers of W. Schmidt \cite{MR778171} and Bombieri and Pila
\cite{MR1016893} for more recent results,  discussions and open questions on this subject.
One may  slightly change Jarn\'{\i}k's framework, and consider the set of integral points which are interpolated by  the graph on $[0,n]$ of  an increasing and strictly convex function satisfying $f(0)=0$ and $f(n)=n$.  It turns out that this question is related to another family of works we shall discuss now. 

In 1979, Arnol'd \cite{arnold_statistics_1980} considered the question of the number of
equivalence classes of convex lattice polygons having a given integer as area
(we say that two polygons having their vertices on $\Z^2$ are equivalent if one
is the image of the other by an automorphism of $\Z^2$).
Later, Vershik changed the constraint in this question and raised the question of
the number, and typical shape, of convex lattice polygons included in a large box
$[-n,n]^2$.
 In 1994, three different solutions to this problem were found by
B\'ar\'any
\cite{barany_limit_1995}, Vershik \cite{MR1275724} and \Sinai{} \cite{MR1283251}. 
Namely, they proved that:
\begin{enumerate}[\quad\bfseries (a)]
    \item The number of convex polygonal chains with vertices in $(\ZZ\cap [0,n])^2$ and joining $(0,0)$ to $(n,n)$ is equal to $\ds \exp[3(\zeta(3)/\zeta(2))^{1/3}n^{2/3}(1+o(1))]$, when $n$ goes to infinity.
    \item The number of vertices constituting a typical line is equivalent, when
$n$ goes to infinity, to $\frac{n^{2/3}}{(\zeta^2(3)\zeta(2))^{1/3}}$.
    \item The limit shape of a typical  convex polygonal line  is the arc of a parabola, which maximizes the affine perimeter.
\end{enumerate}
Note that the approach of \Sinai was recently made rigorous and extended by Bogachev and Zarbaliev \cite{bogachev_zarbaliev_universality_2011}.

Later, Vershik and Zeitouni \cite{MR1679585} proved, for a class of analogous problems, a large deviation principle involving the affine perimeter of the line.
Finally, Acketa and \Zunic, while considering the maximal number of vertices for a lattice polygon included in a square, proved shortly after in \cite{acketa_zunic_maximal_1995} the analog of Jarn\'ik's result, namely that the largest  number of vertices for an increasing convex chain on $\Z_+^2$ of Euclidean length $n$ is asymptotically equivalent to $3\left(\frac{n}{\pi}\right)^{2/3}$.

The nature of these results shows that this problem is related to both 
affine differential geometry and geometry of numbers. Indeed, the parabola
found as limit shape coincides with the convex curve inside the square having the
largest affine perimeter. Furthermore, the appearance of the values of the
Riemann zeta function shows the arithmetic aspects of the problem. One could
show indeed that if the lattice $\Z^2$ was replaced by a Poisson Point Process
having intensity one (which can be thought as the most isotropic ``lattice" one
can imagine), the constants $ (\zeta^2(3)\zeta(2))^{-1/3} \approx 0.749$ and 
$3(\zeta(3)/\zeta(2))^{1/3} \approx 2.702$ would be merely raised respectively to 1 and 3
in probability.

\subsection{Main results}

Our aim in this paper is to improve the three results \textbf{(a),(b),(c)} described above.
In particular, we shall address the following natural extension of \textbf{(c)} which appears as an open question in Vershik's 1994 article:
\begin{quotation}
\itshape
Theorem 3.1 shows how the number of vertices of a typical polygonal line grows. However, one can consider some other fixed growth, say, $\sqrt n$, and look for the limit shapes for uniform distributions connected with this growth [...]
\end{quotation}
One of our results is that, not only there still exists a limit shape when the number of vertices is constrained, but also the parabolic limit shape is actually universal for all growth rates. 
The following theorem is a consequence of Theorem~\ref{thm:limit_shape_numerous} of Section~\ref{sec:limit_shape} and Theorem~\ref{thm:limit_shape_few} of Section~\ref{sec:few}.

\begin{theorem*}
    The Hausdorff distance between a random convex chain on $(\frac{1}{n}\ZZ\cap [0,1])^2$ joining $(0,0)$ to $(1,1)$ with at most $k$ vertices, and the arc of parabola
    \[
        \{(x,y) \in [0,1]^2 : \sqrt{y}+ \sqrt{1-x} = 1\},
    \]
    converges in probability to $0$ when both $n$ and $k$ tend to $+\infty$. 
\end{theorem*}

The proof of this theorem requires a detailed combinatorial analysis of convex chains with a constrained number of vertices. This is the purpose of Theorem~\ref{thm:detailed_comb} in Section~\ref{sec:dca} which generalizes point \textbf{(b)}. We obtain, for any positive number $\c$, a logarithmic equivalent
of the number of lines having roughly $\c\, n^{2/3}$ vertices.
This question is reminiscent of other ones considered, for instance, by \Erdos and Lehner \cite{erdos_lehner_distribution_1941}, Arratia
and Tavar\'e 
\cite{arratia_tavare_independent_1994}, or Vershik and Yakubovich \cite{MR1877604} who were studying
combinatorial objects (integer partitions, permutations, polynomials over finite field, Young
tableaux, etc.) having a specified number of summands 
(according to the setting, we  call summands, cycles, irreducible divisors, etc.).

The method we use emphasizes the connection of the combinatorial analysis with the zeros of the zeta function. We show how the Riemann Hypothesis leads to an asymptotic equivalent of the number of convex chains, improving point \textbf{(a)} above:
\begin{conjecture}
The number $p(n)$ of lattice convex chains in $[0,n]^2$ from $(0,0)$ to $(n,n)$ satisfies
\[
    p(n) \sim \frac{e^{-2\zeta'(-1)}}{(2\pi)^{7/6} \sqrt{3} \kappa^{1/18} n^{17/18}}\exp\Biggl[3\kappa^{1/3} n^{2/3} + \sum_{\substack{\zeta(\rho) = 0\\ \Re(\rho) = \frac{1}{2}}} \frac{\Gamma(\rho)\zeta(\rho+1)\zeta(\rho-1)}{\zeta'(\rho)} \left(\frac{n}{\kappa}\right)^{\rho/3}\Biggr]
\]
where $\kappa = \zeta(3)/\zeta(2)$ and the summation is taken over all zeros of $\zeta$ with real part $\frac{1}{2}$.

\end{conjecture}

\subsection{Organization of the paper}
In Section~\ref{sec:dca}, we detail the  combinatorial aspect of the result of
\cite{barany_limit_1995}, \cite{MR1283251}, \cite{MR1275724} by proving Theorem~\ref{thm:detailed_comb}. Following \Sinai{}'s approach, the method, borrowed from classical
ideas of statistical physics, relies  on the introduction of a grand
canonical ensemble which endows the considered combinatorial object with a parametrized probability measure. Then, the strategy consists in calibrating the parameters of the probability in order to
fit with the constraints one has to deal with. Namely, in  our question, it turns
out that one can add one parameter in \Sinai{}'s probability distribution that makes it
possible to take into account, not only the location of the extreme point of the
chain but also the number of vertices it contains.
In this model, we are able to establish a contour-integral representation of the logarithmic partition function in terms of Riemann's and Barnes' zeta functions. The residue analysis of this representation leads to precise estimates of this function as well as of its derivatives, which correspond to the moments of the random variables of interest such as the position of the terminal point and the number of vertices of the chain. Using a local limit theorem, we finally obtain the
asymptotic behavior of the number of lines having $\c\, n^{2/3}$ vertices in terms of the polylogarithm functions $\Li_1,\Li_2,\Li_3$. We also obtain an asymptotic formula for the number of lines having a number $k$ of vertices satisfying $\log n \ll k \ll n^{2/3}$.

In Section~\ref{sec:limit_shape}, we derive results about the limit shape of lines having a fixed number of vertices $k \gg \log n$, answering the question of Vershik in a wide range.

In Section~\ref{sec:few}, we extend the results about combinatorics and limit shape beyond $\log n$. The approach here is radically different and more elementary, but limited to $k \ll n^{1/3}$. It relies on the comparison with a continuous setting which has been studied by \Barany~\cite{barany_sylvester_1999} and \Barany, Rote, Steiger, Zhang \cite{barany_al_central_2000}.

In Section~\ref{sec:jarnik}, we go back to Jarn\'{\i}k's original problem. In addition to Jarn\'{\i}k's result that we recover, we give the asymptotic number of chains, typical number of vertices, and  limit shape, which is an arc of a circle, in this different framework. 

Furthermore, one may mix both types of conditions and  the statistical physical method still applies. In Section~\ref{sec:onion}, we obtain, for the convex lines joining $(0,0)$ to $(n,n)$ and having  a given total length, a continuous family of convex limit shapes that interpolates the diagonal of the square and the two sides of the square, going through the above arc of parabola and  arc of  circle.

Section~\ref{sec:sf} is devoted to a formal derivation of the above conjecture about the number of convex chains.

\section{A one-to-one  correspondence}
\label{sec:correspondence}

We start this paper by reminding the correspondence between finite convex polygonal
chains issuing from $0$ whose vertices define increasing sequences in both
coordinates and finite distributions of multiplicities on the set of pairs of coprime positive integers.

More precisely,  let $\Pi$ denote the set of finite planar convex polygonal chains $\Gamma$ issuing
from 0 such that the vertices of $\Gamma$ are points of the lattice $\Z^2$ and the angle
between each side of $\Gamma$ and the horizontal axis is in the interval
$[0,\pi/2]$.  Now consider the set $\XX$ of all vectors $x = (x_1,x_2)$ whose coordinates are coprime positive integers including the pairs $(0,1)$ and $(1,0)$. \Sinai{} observed that the space $\Pi$ admits a simple alternative description in terms of distributions of multiplicities on $\XX$.

\begin{lem}[\Sinai{}'s correspondence\cite{MR1283251}]
\label{lem:corresp}
    The space \(\Pi\) is in one-to-one correspondence with the space
    \(\Omega\) of nonnegative integer-valued functions \(x \mapsto
    \omega(x)\) on \(\XX\) with finite support (that is \(\omega(x)\neq
    0\) for only finitely many \(x \in \XX\)).
\end{lem}

The inverse map \(\Omega \to \Pi\) corresponds to the following simple construction: for a given multiplicity distribution \(\omega \in \Omega\) and for all \(\theta\in [0,\infty]\), let us define
\begin{equation}
    X_i^\theta(\omega) := \sum_{\substack{(x_1,x_2) \in \XX\\ x_2 \leq \theta x_1}} \omega(x)\cdot x_i, \qquad i \in \{1,2\}.
\end{equation}
When \(\theta\) ranges over \([0,\infty]\), the function \(\theta \mapsto X^\theta(\omega) = (X^\theta_1(\omega),X^\theta_2(\omega))\) takes a finite number of values which are points of the lattice quadrant \(\ZZ_+^2\). These points are in convex position since we are adding vectors in increasing slope order. The convex polygonal curve \(\Gamma \in \Pi\) associated to \(\omega\) is simply the linear interpolation of these points starting from \( (0,0) \).

\begin{remark*}
    This correspondence is a discrete analogue of the Gauss-Minkowski transformation which was used by Vershik and Zeitouni \cite{MR1679585}.
\end{remark*}

\section{A detailed combinatorial analysis}
\label{sec:dca}

For every \(n = (n_1,n_2) \in \ZZ_+^2\) and \(k\in \ZZ_+\), define \(\Pi(n; k)\) the subset of $\Pi$ consisting of polygonal chains \(\Gamma \in \Pi\) with endpoint \(n\)
and having  \(k\) edges, and denote by \(p(n;k) := \left|\Pi(n;k)\right|\) its cardinality. The restriction of \Sinai{}'s correspondence (see Lemma~\ref{lem:corresp}) to the subspace \(\Pi(n;k)\) induces a bijection with the subset \(\Omega(n;k)\) of \(\Omega\) consisting of multiplicity distributions \(\omega \in \Omega\) such that the "observables"
\[
    X_1(\omega) := \sum_{x \in \XX} \omega(x)\cdot x_1,\quad
    X_2(\omega) := \sum_{x \in \XX} \omega(x)\cdot x_2,\quad 
    K(\omega) := \sum_{x \in \XX} \indicator{\omega(x) > 0}
\]
are respectively equal to \(n_1,n_2\) and \(k\). Notice that $X_1 = X_1^\infty$ and $X_2 = X_2^\infty$ with the notations of the previous section.

Out first theorem gives the asymptotic exponential behavior of $p(n;k)$ in terms of the functions $\c$ and $\e$ defined for all $\lambda \in (0,+\infty)$ by
\begin{gather*}
    \c(\lambda)=\frac{\lambda\LL(1-\lambda)}{1-\lambda}
    \times\frac{1}{\zeta(2)^{1/3}(\zeta(3)- 
    \LLL(1-\lambda))^{2/3}},
    \\
    \e(\lambda)=3\left(\frac{\zeta(3)- 
    \LLL(1-\lambda)}{\zeta(2)}\right)^{1/3}-
    \frac{\lambda\ln(\lambda)\LL(1-\lambda)}{1-\lambda}
    \times\frac{1}{\zeta(2)^{1/3}(\zeta(3)- 
    \LLL(1-\lambda))^{2/3}}.
\end{gather*}

\begin{figure}[h]
    \begin{center}
        \includegraphics[scale=1.2]{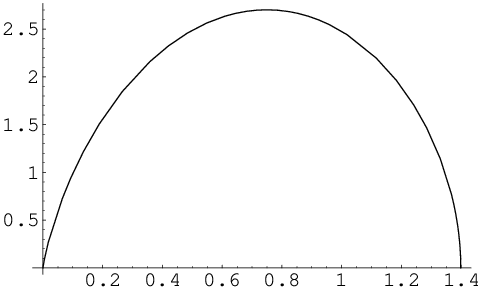}
        \caption{Graph of $(\c(\lambda),\e(\lambda))$. The point of maximal $\e$-coordinate corresponds to typical chains. The point of maximal $\c$-coordinate correspond to chains with a maximal number of vertices.}
    \end{center}
\end{figure}

\begin{theorem}
    \label{thm:detailed_comb}
    Suppose that $|n|$ and $k$ tend to $+\infty$ such that $n_1 \asymp n_2$ and $\log |n|$ is asymptotically negligible compared to $k$.
    
    \begin{itemize}
        \item 
            If there exists $\lambda \in (0,+\infty)$ such that $k \sim \c(\lambda)(n_1n_2)^{1/3}$, then
    \[
        \log p(n;k) \sim  \e(\lambda) (n_1n_2)^{1/3}.
    \]
        \item  If $k$ is asymptotically negligible compared to $(n_1n_2)^{1/3}$, then
            \[
                p(n;k) = \left(\frac{n_1n_2}{k^3}\right)^{k+o(k)}.
            \]
    \end{itemize}
\end{theorem}

\begin{remark}
    It will appear in the core of the proof that one cannot obtain additional terms in the expansion of $\log p(n;k)$ without strong knowledge of the localization of the zeros of Riemann'z zeta function. We show in the last section of the paper how the result can be improved up to an asymptotic equivalent under the assumption of the Riemann Hypothesis.
\end{remark}

\begin{remark}
    The function $\e(\lambda)$ is maximal for $\lambda=1$, that is to say when there is no penalization. The corresponding coefficients are
    \begin{eqnarray*}
        \c(1)=\dfrac{1}{(\zeta(2)\zeta(3)^2)^{1/3} },
         &
    \e(1)=3\left(\dfrac{\zeta(3)}{\zeta(2)}\right)^{1/3},
    \end{eqnarray*}
    which recovers the results of \cites{MR1283251,MR1275724,barany_limit_1995}. 
\end{remark}

\begin{remark}
	\label{rem:longest}
	As a byproduct of Theorem~\ref{thm:detailed_comb}, one can deduce the asymptotic behavior of the maximal number of integral points that an increasing convex function satisfying $f(0)=0$ and $f(n)=n$ can interpolate. This question and its counterpart, concerning the maximal convex lattice polygons inscribed in a convex set was solved by Acketa and \Zunic \cite{acketa_zunic_maximal_1995} who proved that $M(n) := \max_{\Gamma \in\Pi(n,n)} \bigl | \Gamma \cap \ZZ^2 \bigr | \sim
3\left(\pi^{-1}n\right)^{2/3}$.

Starting from Theorem~\ref{thm:detailed_comb}, the proof goes as follows. We first notice that $\e(\lambda)$ vanishes when $\lambda$ goes to
infinity. In the same time,
$
    \c(\lambda) {\sim}
    \frac{-\LL(1-\lambda)}{\zeta(2)^{1/3}(- 
\LLL(1-\lambda))^{2/3}}
$
which tends to ${3}\pi^{-2/3}$.
Since $\e(\lambda)$ remain strictly positive, we get $\liminf M(n) n^{-2/3}\geq
    \frac{3}{\pi^{2/3}}$
    Now, let $\epsilon>0$, and suppose 
$
    \limsup M(n) n^{-2/3}
    \geq \frac{3(1+2\epsilon)}{\pi^{2/3}}.
$
Then, for arbitrary large $n$, there is a chain $\Gamma\in\Pi(n,n)$ having cardinality $[(1+\epsilon)3(n/\pi)^{2/3}]$. By choosing $k = 3(\pi^{-1}n)^{2/3}$ vertices on this chain, we get already a subset of $\Pi(n; k)$ whose cardinality is $e^{cn^{2/3}}$ with $c > 0$. This enters in contradiction with the fact that $\lim_{\lambda\to\infty} \e(\lambda) = 0$.
\end{remark}

\subsection{Modification of \Sinai{}'s model and proof of Theorem~\ref{thm:detailed_comb}}

Let \(\L = \{(x_1,x_2) \in \Z_+^2 : \gcd(x_1,x_2) = 1\}\). For any \(\omega \in \Omega := \mathcal{F}(\L,\Z_+)\), let \(K(\omega)\) be the number of \(x\in\L\) such that \(\omega(x) > 0\), that is to say
\[
    K(\omega) := \sum_{x\in\L} \indicator{\omega(x) > 0}.
\]

For all $\lambda > 0$ and for every couple of parameters $\beta = (\beta_1,\beta_2) \in (0,+\infty)^2$, we endow \(\Omega\) with the Gibbs measure defined for $\omega \in \Omega$ by
\begin{align*}
    \PPbl(\omega) & := \frac{1}{\Zbl}\exp\left[- \sum_{x \in \L} \omega(x)\, \beta \cdot x\right] \lambda^{K(\omega)}\\
                  & = \frac{1}{\Zbl} e^{-\beta_1 X_1(\omega)} e^{-\beta_2 X_2(\omega)} \lambda^{K(\omega)},
\end{align*}
where the \emph{partition function} $\Zbl$ is chosen as the normalization constant
\begin{equation}
	\label{eq:Zbl}
    \Zbl = \sum_{n \in \Z_+^2} \sum_{k\geq1} p(n ; k)\, e^{-\beta \cdot n} \lambda^k.
\end{equation}
Note that $\Zbl$ is finite for all values of the parameters $(\beta,\lambda) \in (0,+\infty)^3$. Indeed, if we denote by $p(n) = \sum_{k\geq 1} p(n;k)$ the total number of convex chains of $\Pi$ with end point $n = (n_1,n_2)$ and $M_n$ the maximal number of edges of such a chain, the following bound holds:
\[
    \Zbl \leq \sum_{n \in \Z_+^2} p(n)\,  \max(1,\lambda)^{M_n} \, e^{-\beta \cdot n}.
\]
We use now the results of \cites{barany_limit_1995,MR1283251,MR1275724} according to which $\log p(n)  = O(|n|^{2/3})$ and of \cite{acketa_zunic_maximal_1995} where Acketa and \Zunic have proven that $M_n = O(|n|^{2/3})$.
We will use in the sequel the additional remark that $\Zbl$ is an analytic function of $\lambda$ for all $\beta > 0$.

Taking $\lambda=1$, the probability $\PPbl$ is nothing but the two-parameter probability distribution introduced by \Sinai{}~\cite{MR1283251}. Under the measure $\PPbl$, the variables $(\omega(x))_{x\in \L}$
are still independent, as in \Sinai{}'s framework,  but follow a geometric distribution only for
$\lambda=1$. In the general case, the measure \(\PPbl\) is absolutely continuous with respect to \Sinai{}'s measure with density proportional to $\lambda^{K(\cdot)}$ and the distribution of $\omega(x)$ is a biased geometric distribution. Loosely speaking, $\PPbl$ corresponds to the introduction of a {\it penalization} 
of the probability by a factor $\lambda$ each time a vertex appears.

Since $\PPbl(\omega)$ depends only on the values of $X_1(\omega)$,
$X_2(\omega)$, and $K(\omega)$, we deduce
that the conditional distribution it induces on $\Omega(n_1,
n_2; k)$ is uniform. For instance, we have the following formula for all $(\beta,\lambda) \in (0,+\infty)^2\times (0,+\infty)$ which will be instrumental in the proof:
\[
    p(n_1,n_2;k)=Z(\beta,\lambda)\, e^{\beta_1 n_1} e^{\beta_2 n_2} \lambda^{-k} 
 \,\PPbl[X_1 = n_1, X_2 = n_2, K =  k].
\]

In order to get a logarithmic equivalent of $ p(n_1,n_2;k)$, our strategy is
to choose the three parameters so that 
\[
    \EEbl\left[X_1\right] = n_1,\quad
    \EEbl\left[X_2\right] = n_2,\quad
    \EEbl\left[K\right]=k.
\]
This will indeed lead to an asymptotic equivalent of $\PPbl[X_1 = n_1, X_2 = n_2, K =  k]$ due to a local limit result. This equivalent having polynomial decay, it will not interfere with the estimation of $\log p(n_1,n_2;k)$. Together with the analysis of the partition function, this local limit result will constitute the heart of the proof.

\subsection{Estimates of the logarithmic partition function and its derivatives}

We need in the following, the analogue to the Barnes bivariate zeta function defined for $\beta=(\beta_1,\beta_2) \in (0,+\infty)^2$ by
\[
    \zeta_2^*(s ; \beta) := \sum_{x \in\L} (\beta_1 x_1 + \beta_2 x_2)^{-s},
\]
this series being convergent for \(\Re(s) > 2\). The following preliminary lemma gives useful properties of this function. This will be done by expressing this function in terms of the Barnes zeta function $\zeta_2(s,w ; \beta)$ which is defined by analytic continuation of the series
$$
\zeta_2(s,w ; \beta) = \sum_{n \in \Z_+^2} (w + \beta_1 n_1 + \beta_2 n_2)^{-s}, \qquad \Re(s) > 2, \Re(w) > 0.
$$
It is well known that $\zeta_2(s,w ; \beta)$ has a meromorphic continuation to the complex $s$-plane with simple poles at $s = 1$ and $2$, and that the residue at $s = 2$ is simply \((\beta_1\beta_2)^{-1}\). In the next lemma, we derive the relation between $\zeta_2$ and $\zeta_2^*$, and we also establish an explicit meromorphic continuation of $\zeta_2$ to the half-plane $\Re(s) > 1$ in order to obtain later polynomial bounds for $|\zeta_2^*(s)|$ as $|\Im(s)| \to +\infty$. 

\begin{lem}
\label{lem:meromorphic_continuation}
The functions $\zeta_2(s,w ; \beta)$ and $\zeta_2^*(s ; \beta)$ have a meromorphic continuation to the complex plane.
\begin{enumerate}[(i)]
    \item The meromorphic continuation of $\zeta_2(s,w ; \beta)$ to the half-plane $\Re(s) > 1$ is given by
\begin{align*}
\zeta_2(s,w ;\beta) &= \frac{1}{\beta_1\beta_2}\frac{w^{-s+2}}{(s-1)(s-2)} + \frac{(\beta_1 + \beta_2)w^{-s+1}}{2\beta_1\beta_2(s-1)} + \frac{w^{-s}}{4}\\
& - \frac{\beta_2}{\beta_1} \int_0^{+\infty} \frac{\{y\}-\frac{1}{2}}{(w+\beta_2 y)^s}\,dy - \frac{\beta_1}{\beta_2} \int_0^{+\infty} \frac{\{x\}-\frac{1}{2}}{(w+\beta_1 x)^s}\,dx\\
& - s\frac{\beta_2}{2} \int_0^{+\infty} \frac{\{y\}-\frac{1}{2}}{(w+\beta_2 y)^{s+1}}dy
- s\frac{\beta_1}{2} \int_0^{+\infty} \frac{\{x\}-\frac{1}{2}}{(w+\beta_1 x)^{s+1}}dx\\
& + s(s+1)\beta_1\beta_2 \int_0^{+\infty}\int_0^{+\infty} \frac{(\{x\} - \frac{1}{2})(\{y\}-\frac{1}{2})}{(w+\beta_1 x +\beta_2 y)^{s+2}}\,dxdy.
\end{align*}
\item The meromorphic continuation of $\zeta_2^*(s;\beta)$ is given for all $s \in \CC$ by
$$
\zeta_2^*(s;\beta) =  \frac{1}{\beta_1^s} + \frac{1}{\beta_2^s} + \frac{\zeta_2(s,\beta_1+\beta_2 ; \beta)}{\zeta(s)}.
$$
\end{enumerate}
\end{lem}

\begin{proof}[Proof of (i)]
Let $\{x\} = x - \lfloor x \rfloor$ denote the fractional part of $x$. We apply the Euler-Maclaurin formula to the partial summation defined by $F(x) = \sum_{n_2\geq0} (w +\beta_1 x + \beta_2 n_2)^{-s}$, leading to
$$
\sum_{n_1 \geq 1} F(n_1) = \int_0^{\infty} F(x)\,dx - \frac{F(0)}{2} + \int_0^\infty (\{x\}-\frac{1}{2})F'(x)\,dx.
$$
We use again the Euler-Maclaurin formula for each of the summations in $n_2$. 
\end{proof}

\begin{proof}[Proof of (ii)]
    Let us express $\zeta_2^*(s;\beta)$ in terms of $\zeta_2(s,\beta_1+\beta_2;\beta)$ for all $s$ with real part $\Re(s) > 2$. The result will follow from the analytic continuation principle. By definition of $\zeta_2^*(s;\beta)$,
$$
\zeta_2^*(s;\beta) - \frac{1}{\beta_1^s} - \frac{1}{\beta_2^s} = \sum_{x_1,x_2 \geq 1} \frac{1}{(\beta_1 x_1 + \beta_2 x_2)^s} 1_{\gcd(x_1,x_2) = 1}.
$$
Using the classical \Moebius{} function $\mu(d)$ taking values in $\{-1,0,1\}$ and the \Moebius{} inversion formula (see \cite{MR2445243}), we obtain 
$$
1_{\gcd(x_1,x_2) = 1} = \sum_{d = 1}^\infty \mu(d) 1_{d \mid x_1} 1_{d \mid x_2},
$$
so we can write the latter expression as
$$
\sum_{d\geq 1} \sum_{x_1,x_2 \geq 1} \frac{\mu(d)}{(\beta_1 x_1 + \beta_2 x_2)^s} 1_{d \mid x_1} 1_{d \mid x_2}= \sum_{d \geq 1} \frac{\mu(d)}{d^s} \sum_{x'_1,x'_2 \geq 1} \frac{1}{(\beta_1 x_1 + \beta_2 x_2)^s}.
$$
Finally, the classical formula
$$
\sum_{d=1}^\infty \frac{\mu(d)}{d^s} = \frac{1}{\zeta(s)}
$$
which holds for all $s$ having real part $\Re(s) > 1$ implies (ii).
\end{proof}

Now we make the connection between these zeta functions and the logarithmic partition function of our modified \Sinai's model. Let us recall that the polylogarithm function $\Li_s(z)$, also known as Jonqui\`ere function, is defined for all complex number $s \in \CC$ by analytic continuation of the series
\[
    \Li_s(z) = \sum_{k=1}^\infty \frac{z^k}{k^s}, \qquad |z| < 1.
\]
For our purpose, the continuation given by the Bose-Einstein integral
\[
    \Li_s(z) = \frac{1}{\Gamma(s)}\int_0^\infty \frac{z t^{s-1}}{e^t - z} dt
\]
for $\Re(s) > 0$ and $z \in \CC\setminus[1,+\infty)$ will be sufficient.

\begin{lem}
\label{lem:integral_partition_function}
Let $c > 2$. For all parameters $(\beta,\lambda) \in (0,+\infty)^2 \times (0,+\infty)$,
$$
\log \Zbl = \frac{1}{2i\pi} \int_{c - i\infty}^{c + i\infty} (\zeta(s+1) - \Li_{s+1}(1-\lambda))\zeta_2^*(s;\beta)\Gamma(s)ds.
$$
\end{lem}

\begin{proof}
Given the product form of the distribution $\PPbl$, we see that the random variables $\omega(x)$ for $x \in\L$ are mutually independent. Moreover, the marginal distribution of $\omega(x)$ is a biased geometric distribution. It is absolutely continuous with respect to the geometric distribution of parameter $e^{-\beta \cdot x}$ with density proportional to $k \mapsto \lambda^{1_{k > 0}}$. In other words, for all $k \in \Z_+$,
$$
\PPbl[\omega(x) = k] = Z_x(\beta, \lambda)^{-1} e^{-k \beta\cdot x} \lambda^{1_{k > 0}}
$$
where the normalization constant $Z_x(\beta, \lambda) = 1 + \lambda\dfrac{e^{-\beta \cdot x}}{1-e^{-\beta \cdot x}}$ is easily computed. We can now deduce the following product formula for the partition function:
\[
    \Zbl = \prod_{x\in\L} Z_x(\beta, \lambda) = \prod_{x\in\L} \left(1 + \lambda\frac{e^{-\beta \cdot x}}{1-e^{-\beta \cdot x}}\right).
\]
For now, we assume that $\lambda \in (0,1)$. Taking the logarithm of the product above
\begin{align*}
    \log \Zbl & = \sum_{x \in \Lstar} \log \left(1 + \lambda\frac{e^{-\beta \cdot x}}{1-e^{-\beta \cdot x}}\right) \\
    &= \sum_{x\in \Lstar} \log(1-(1-\lambda)e^{-\beta\cdot x}) - \sum_{x\in\Lstar} \log(1-e^{-\beta\cdot x}) \\
    &= \sum_{x\in\Lstar} \sum_{r \geq 1} \frac{1 - (1-\lambda)^r}{r} e^{-r\beta \cdot x}.
\end{align*}
Now we use the fact that the Euler gamma function $\Gamma(s)$ and the exponential function are related through Mellin's inversion formula
$$
e^{-z} = \frac{1}{2i\pi} \int_{c - i\infty}^{c + i\infty} \Gamma(s) z^{-s} ds,
$$
for all $c > 0$ and $z \in \CC$ with positive real part. Choosing $c > 2$ so that the series and the integral all converge and applying the Fubini theorem, this yields
\begin{align*}
    \log \Zbl & = \frac{1}{2i\pi} \sum_{x\in\Lstar} \sum_{r\geq 1} \int_{c - i\infty}^{c + i\infty}\frac{1-(1-\lambda)^r}{r} r^{-s} (\beta \cdot x)^{-s} \Gamma(s) ds \\
    & = \frac{1}{2i\pi} \int_{c - i\infty}^{c + i\infty} (\zeta(s+1) - \Li_{s+1}(1-\lambda))\zeta_2^*(s;\beta)\Gamma(s)\,ds.
    \end{align*}

    The lemma is proven for all $\lambda \in (0,1)$. The extension to $\lambda > 0$ will now result from analytic continuation. We already noticed that the left hand term is analytic in $\lambda$ for all fixed $\beta$. Proving the analyticity of the right hand term requires only to justify the absolute convergence of the integral on the vertical line. From Lemma~\ref{lem:meromorphic_continuation}, we know that $\zeta_2^*(c + i\tau ; \beta)$ is polynomially bounded as $|\tau|$ tends to infinity. Taking $s = c - 1 + i\tau$, successive integrations by parts of the formula
    \[
    	(\zeta(s+1) - \Li_{s+1} (1-\lambda))\Gamma(s+1) = \lambda \int_0^\infty \frac{e^x x^s}{(e^x - 1)(e^x - 1 + \lambda)} \,dx
    \]
    show for all integer $N > 0$, there exists a constant $C_N > 0$ such that, uniformly in $\tau$, 
    \begin{equation}
    	\label{eq:decay_lambda}
    	\bigl | (\zeta(s+1) - \Li_{s+1} (1-\lambda))\Gamma(s+1) \bigr | \leq \frac{ C_N\lambda}{(1+|\tau|)^N}.
    \end{equation}
\end{proof}

Finally, the next Lemma makes use of the contour integral representation of $\log Z(\beta,\lambda)$ to derive at the same time an asymptotic formula for each one of its derivatives.

\begin{lem}
\label{lem:derivatives_Z}
Let $(p,q_1,q_2) \in \ZZ_+^3$. For all $\epsilon > 0$, there exists $C > 0$ such that
$$
\left|\left[\lambda\frac{\partial}{\partial \lambda}\right]^p\left[\frac{\partial}{\partial \beta_1}\right]^{q_1}\left[\frac{\partial}{\partial \beta_2}\right]^{q_2} \left( \log \Zbl - \frac{\zeta(3)-\Li_3(1-\lambda)}{\zeta(2)\beta_1\beta_2}\right)\right| \leq \frac{C\,\lambda}{|\beta|^\kappa}
$$
with $\kappa = q_1 + q_2 + 1 + \epsilon$, uniformly in the region $\{(\beta,\lambda) : \epsilon < \frac{\beta_1}{\beta_2} < \frac{1}{\epsilon} \text{ and } 0 < \lambda < \frac{1}{\epsilon}\}$. 
\end{lem}

\begin{proof}
    Lemma~\ref{lem:integral_partition_function} provides an integral representation of the logarithmic partition function $\log\Zbl$. We will use the residue theorem to shift the contour of integration from the vertical line $\Re(s) = 3$ to the line $\Re(s) = 1 + \epsilon$. Lemma~\ref{lem:meromorphic_continuation} shows that the function $M(s) := (\zeta(s+1)-\Li_{s+1}(1-\lambda)))\zeta_2^*(s;\beta)\Gamma(s)$ is meromorphic in the strip $1 < \Re(s) < 3$ with a single pole at $s=2$, where the residue is given by
$$
\frac{\zeta(3)-\Li_3(1-\lambda)}{\zeta(2)} \cdot \frac{1}{\beta_1\beta_2}
$$ 
From the inequality~\eqref{eq:decay_lambda}, Lemma~\ref{lem:meromorphic_continuation} and the fact that $|\zeta(s)|$ has no zero with $\Re(s) > 1$, we see that $M(s)$ vanishes uniformly in $1 + \epsilon \leq \Re(s) \leq 3$ when $|\Im(s)|$ tends to $+ \infty$. By the residue theorem,
\begin{equation}
\label{eq:residue_theorem}
\log \Zbl= \frac{\zeta(3)-\Li_3(1-\lambda))}{\zeta(2)\beta_1\beta_2}+ \frac{1}{2i\pi} \int_{1 + \epsilon - i\infty}^{1 + \epsilon + i\infty} M(s)\,ds.
\end{equation}
From the Leibniz rule applied in the formula of Lemma~\ref{lem:meromorphic_continuation}~(i), we obtain directly the meromorphic continuation of $\frac{\partial^{q_1}}{\partial \beta_1^{q_1}}\frac{\partial^{q_2}}{\partial \beta_2^{q_2}}\zeta_2(s,\beta_1+\beta_2;\beta)$ in the half-plane $\Re(s) > 1$. We also obtain the existence of a constant $C > 0$ such that
\[
    \left|\left[\frac{\partial}{\partial \beta_1}\right]^{q_1}\left[\frac{\partial}{\partial \beta_2}\right]^{q_2}\zeta_2(1 + \epsilon + i\tau ,\beta_1+\beta_2;\beta)\right| \leq  \frac{C\, |\tau|^{2 + q_1 + q_2}}{|\beta|^\kappa}
\]
with $\kappa = q_1 + q_2 + 1 + \epsilon$.
A reasoning similar to the one we have used in order to derive \eqref{eq:decay_lambda} shows that for all integers $p$ and $N > 0$, there exists a constant $C_{p,N}$ such that, uniformly in $\tau$, 
    \[
    	\left | \left[\lambda \frac{\partial}{\partial\lambda}\right]^p (\zeta(s+1) - \Li_{s+1} (1-\lambda))\Gamma(s+1) \right | \leq \frac{ C_{p,N}\,\lambda}{(1+|\tau|)^N}.
    \]

    In order to differentiate both sides of equation~\eqref{eq:residue_theorem} and permute the partial derivatives and the integral sign, we have to mention the fact that the Riemann zeta function is bounded from below on the line $\Re(s) = 1 + \epsilon$ and that the derivatives of $\Li_s(1-\lambda)$ with respect to $\lambda$ are all bounded. This also gives the announced bound on the error term.
\end{proof}

\begin{remark}
    As shown by this proof, the asymptotic expansion of $\log Z$ is directly related to the localization of the zeros of Riemann's zeta function in the critical strip $0 < \Re(s) \leq 1$. For instance, using the fact that there is no zero on the line $\Re(s) = 1$ (which is a form of the Prime Number Theorem), one may improve the exponent $\kappa$ in the remainder to $\kappa = q_1 + q_2 + 1$. On the other hand, the existence of infinitely many zeros on the line $\Re(s) = \frac{1}{2}$ implies that $\log Z$ presents fluctuations of order at least $|\beta|^{-1/2}$. Finally, finding the next term in the expansion depends on the Riemann Hypothesis.
\end{remark}

\subsection{Calibration of the shape parameters}
\label{sec:calibration}

When governed by the Gibbs measure \(\PPbl\), the expected value of the random vector with components
\[
    X_1(\omega) = \sum_{x\in\L} \omega(x) x_1, \quad
    X_2(\omega) = \sum_{x\in\L} \omega(x) x_2, \quad
    K(\omega) = \sum_{x\in\L} \indicator{\omega(x) > 0},
\]
is simply given by the logarithmic derivatives of the partition function \(\Zbl\). Remember that we planned to choose \( \lambda \) and  \( \beta_1, \beta_2 \) as functions of \(n = (n_1,n_2)\) an \(k\) in order for the probability $\PP[X_1 = n_1, X_2 = n_2, K = k]$ to be maximal, which is equivalent to $\EE(X_1)=n_1$, $\EE(X_2)=N_2$ and $\EE(K)=k$. We address this question in the next lemma.

\begin{lem}
    \label{lem:parameters}
    Assume that $n_1,n_2,k$ tend to infinity with $n_1 \asymp n_2$ and $|k| = O(|n|^{2/3})$. There exists a unique choice of $(\beta_1,\beta_2,\lambda)$ as functions of $(n,k)$ such that
    \[
        \EEbl[X_1] = n_1, \quad \EEbl[X_2] = n_2, \quad \EEbl[K] = k.
    \]
    Moreover, they satisfy
    \begin{equation}
        \label{eq:parameters}
        n_1 \sim \frac{\zeta(3) - \Li_3(1-\lambda)}{\zeta(2)({\beta_1})^2{\beta_2}}, \quad n_2 \sim \frac{\zeta(3)-\Li_3(1-\lambda)}{\zeta(2){\beta_1}({\beta_2})^2}, \quad
        k \sim -\frac{\lambda\partial_\lambda
             \Li_3(1-\lambda)}{\zeta(2){\beta_1}{\beta_2}}.
     \end{equation}
    If $k = o(|n|^{2/3})$, then $\lambda$ goes to $0$ and the above relations yield
    \[
        \beta_1 \sim \frac{k}{n_1}, \quad \beta_2 \sim \frac{k}{n_2}, \quad \lambda \sim \frac{k^3}{n_1n_2}.
    \]
\end{lem}

\begin{proof}
    With the change of variable $\lambda = e^{-\gamma}$, the existence and uniqueness of $(\beta,\lambda)$ are equivalent to the fact that the function
\[
f \colon (\beta_1,\beta_2,\gamma) \mapsto \beta_1 n_1 + \beta_2 n_2 + \gamma k + \log Z(\beta, e^{-\gamma})
\]
has a unique critical point in the open domain $D = (0,+\infty)^2 \times \RR$.
First observe that $f$ is smooth and strictly convex since its Hessian matrix is actually the covariance matrix of the random vector $(X_1,X_2,K)$. In addition, from the very definition~\eqref{eq:Zbl} of $\Zbl$, we can see that $f$ converges to $+\infty$ in the neighborhood of any point of the boundary of $D$ as well as when $|\beta_1| + |\beta_2| + |\gamma|$ tends to $+\infty$. The function being continuous in $D$, this implies the existence of a minimum, which by convexity is the unique critical point $(\beta^*,\gamma^*)$ of $f$.

From now on, we will be concerned and check along the proof that we stay in the regime $\beta_1,\beta_2\to 0$, $\beta_1 \asymp \beta_2$, and $\gamma$ bounded from below. From Lemma~\ref{lem:derivatives_Z}, we can approximate $f$ by the simpler function
\[
    g \colon (\beta_1,\beta_2,\gamma) \mapsto \beta_1 n_1 + \beta_2 n_2 + \gamma k + \dfrac{\zeta(3)-\Li_3(1-e^{-\gamma})}{\beta_1\beta_2}
\]
with $|f(\beta,\gamma) - g(\beta,\gamma)| \leq \dfrac{Ce^{-\gamma}}{|\beta|^{3/2}}$ for some constant $C > 0$. The unique critical point $(\tilde{\beta},\tilde{\gamma})$ of $g$ satisfies
\[
        n_1 = \frac{\zeta(3) - \Li_3(1-e^{-\tilde{\gamma}})}{\zeta(2)(\tilde{\beta_1})^2\tilde{\beta_2}}, \quad n_2 = \frac{\zeta(3)-\Li_3(1-e^{-\tilde{\gamma}})}{\zeta(2)\tilde{\beta_1}(\tilde{\beta_2})^2}, \quad k = -\frac{e^{-\tilde{\gamma}}\partial_\lambda
         \Li_3(1-e^{-\tilde{\gamma}})}{\zeta(2)\tilde{\beta_1}\tilde{\beta_2}}.
\]
The aim now is to prove that $(\beta^*,\gamma^*)$ is close to $(\tilde{\beta},\tilde{\gamma})$. To this aim, we find a convex neighborhood $C$ of $(\tilde{\beta},\tilde{\gamma})$ such that $g|_{\partial C} \geq g(\tilde{\beta},\tilde{\gamma}) + \frac{Ce^{-\tilde{\gamma}}}{\tilde{\beta_1}\tilde{\beta_2}}$ .
In the neighborhood of $(\tilde{\beta},\tilde{\gamma})$ the expression of the Hessian matrix of $g$ yields $g(\tilde{\beta}_1 + t_1, \tilde{\beta}_2 + t_2, \tilde{\gamma} + u) \geq g(\tilde{\beta}_1 , \tilde{\beta}_2 , \tilde{\gamma} ) + \frac{\tilde{C}e^{-\tilde{\gamma}} }{(\tilde{\beta}_1\tilde{\beta}_2)^2}(\|t\|^2 + \tilde{\beta}_1\tilde{\beta}_2 |u|^2)$. Therefore we need only take
\[
    C = [\tilde{\beta}_1 - C_1\tilde{\beta}_1^{5/4}, \tilde{\beta}_1 + C_1\tilde{\beta}_1^{5/4}] \times  
    [\tilde{\beta}_2 - C_2\tilde{\beta}_2^{5/4}, \tilde{\beta}_2 + C_2\tilde{\beta}_2^{5/4}] \times [\tilde{\gamma} - C_3 |\beta|^{1/4}, \tilde{\gamma} + C_3 |\beta|^{1/4}].
\]
Therefore, $f|_{\partial C} > f(\tilde{\beta},\tilde{\gamma})$. By convexity of $f$ and $C$ this implies $(\beta^*,\gamma^*) \in C$. Hence
\[
    \beta_1^* \sim \tilde{\beta}_1, \quad
    \beta_2^* \sim \tilde{\beta}_2, \quad
    e^{-\gamma^*} \sim e^{-\tilde{\gamma}},
\]
concluding the proof.
\end{proof}

\subsection{A local limit theorem}

In this section, we show that the random vector $(X_1,X_2,K)$ satisfies a
local limit theorem when the parameters are calibrated as above. Let $\Gammabl$ be the covariance matrix under the measure $\PPbl$ of the random vector $(X_1,X_2, K)$.
\begin{theorem}[Local limit theorem]
    \label{thm:local_limit}
    Let us assume that $n_1,n_2,k$ tend to infinity such that $n_1 \asymp n_2 \asymp |n|$, $\log |n| = o(k)$, and $k = O(|n|^{2/3})$. For the choice of parameters made in Lemma~\ref{lem:parameters},
    \begin{equation}
        \label{eq:local_limit}
        \PPbl[X = n, K = k] \sim \frac{1}{(2\pi)^{3/2}}\frac{1}{\sqrt{\det \Gammabl}}.
    \end{equation}
Moreover,
\begin{equation}
    \det \Gammabl \asymp \frac{|n|^4}{k}
\end{equation}
If $k = o(|n|^{2/3})$,
    \begin{equation}
        \PPbl[X = n, K = k] \sim \frac{1}{(2\pi)^{3/2}}\frac{\sqrt{k}}{n_1n_2}
    \end{equation}
\end{theorem}
 
This result is actually an application of a more general lemma proven by the first author in \cite{bureaux_partitions_2014}*{Proposition 7.1}. In order to state the lemma, we introduce some notations. Let $\sigmabl^2$ be the smallest eigenvalue of $\Gammabl$.  
Introducing $X_{1,x} = \omega(x) \cdot x_1$, $X_{2,x} = \omega(x) \cdot x_2$ and $K_x = 1_{\{\omega(x) > 0\}}$ as well as $\overline{X_{1,x}} , \overline{X_{2,x}},\overline{K_x}$ their centered counterparts, let $\Lbl$ be the Lyapunov ratio 
\[
\Lbl := \sup_{(t_1,t_2,u)\in \RR^3}  \sum_{x \in \L} \frac{\EEbl \left\lvert t_1 \overline{X_{1,x}} + t_2 \overline{X_{2,x}} +  u\overline{K_x}\right\rvert^3}{\Gammabl(t_1,t_2,u)^{3/2}}.
\]where $\Gammabl(\cdot)$ stands for the quadratic form canonically associated to $\Gammabl$.
Let $ \phibl(t,u) = \EEbl (e^{i(t_1 X_1 + t_2 X_2 + uK})$ for all $(t_1,t_2,u) \in \RR^3$. Finally, we consider the ellipsoid $\Ecal_{\beta,\lambda}$ defined by
\[
\mathcal{E}_{\beta,\lambda} := \left\{(t_1,t_2,u) \in \RR^3 : \Gammabl(t_1,t_2,u) \leq (4\Lbl)^{-2}\right\}.
\]

The following lemma is a reformulation of Proposition~{7.1} in \cite{bureaux_partitions_2014}. It gives three conditions on the product distributions $\PPbl$ that entail a local limit theorem with given speed of convergence.

 \begin{lem}
     \label{lem:framework}
     With the notations introduced above, suppose that there exists a family of number $(a_{\beta,\lambda})$ such that
     \begin{gather}
     \frac{1}{\sigmabl\sqrt{\det \Gammabl}} = \bigO(a_{\beta,\lambda}), \\
              \frac{\Lbl}{\sqrt{\det \Gammabl}} = \bigO(a_{\beta,\lambda}), \\
              \label{eq:cramer}
               \sup_{(t,u) \in [-\pi,\pi]^3\setminus \Ecal_{\beta,\lambda}} \left|\phibl(t,u)\right| = \bigO(a_{\beta,\lambda}).
     \end{gather}
     Then, a local limit theorem holds uniformly for $\PPbl$ with rate $a_{\beta,\lambda}$:
     \[
         \sup_{(n,k)\in \Z^3}\; \left|\PPbl[X = n, K = k] - \frac{\exp\left[-\frac{1}{2}\Gammabl^{-1}\bigl((n,k) - \EEbl (X,K)\bigr)\right]}{(2\pi)^{3/2}\sqrt{\det \Gammabl}}\right| = \bigO(a_{\beta,\lambda}).
     \]
 \end{lem}
 
 When governed by the Gibbs measure \(\PPbl\), the covariance matrix $\Gammabl$ of the random vector \((X_1,X_2,K)\) is simply given by the Hessian matrix of the log partition function $\log \Zbl$. Let \(u(\lambda) := (\zeta(3) - \Li_3(1-\lambda))/\zeta(2)\) for \(\lambda > 0\). Applications of Lemma~\ref{lem:derivatives_Z} for all $(p,q_1,q_2) \in \ZZ_+^3$ such that $p+q_1+q_2 = 2$ imply that this covariance matrix is asymptotically equivalent to
\[
    \begin{bmatrix}
        \beta_1\beta_2 & 0 & 0\\
        0 & \beta_1^3 \beta_2 & 0\\
        0 & 0 & \beta_1\beta_2^3
    \end{bmatrix}^{-\frac{1}{2}}
    \begin{bmatrix}
        \lambda^2 u''(\lambda) + \lambda u'(\lambda)  &
        \lambda u'(\lambda) &
        \lambda u'(\lambda) \\

        \lambda u'(\lambda) &
        2u(\lambda) &
        u(\lambda) \\

        \lambda u'(\lambda) &
        u(\lambda) &
        2u(\lambda) &
    \end{bmatrix}
    \begin{bmatrix}
        \beta_1\beta_2 & 0 & 0\\
        0 & \beta_1^3 \beta_2 & 0\\
        0 & 0 & \beta_1\beta_2^3
    \end{bmatrix}^{-\frac{1}{2}}.
\]
A straightforward calculation shows that this matrix is positive definite for all $\lambda > 0$.

\begin{lem}
\label{lem:covariance}
The random vector \((X_1,X_2, K)\) has a covariance matrix \(\Gammabl\) satisfying
\[
\Gammabl(t,u) \asymp \frac{(n_1)^{5/3}}{(\lambda n_2)^{1/3}} |t_1|^2 + \frac{(n_2)^{5/3}}{(\lambda n_1)^{1/3}} |t_2|^2 + (\lambda n_1n_2)^{1/3} |u|^2, \qquad |n| \to +\infty.
\]
\end{lem}

\begin{proof}
All the coefficients of the previous matrix $u(\lambda), \lambda u'(\lambda), \lambda^2 u''(\lambda)$ are of order $\lambda$ in the neighborhood of $0$, and the determinant is equivalent to $\lambda^3$. Therefore, the eigenvalues are also of order $\lambda$. The result follows from the fact that the values of $\beta_1$ and $\beta_2$ are given by~\eqref{eq:parameters} and that $\zeta(3)-\Li_3(1-\lambda) \asymp \zeta(2) \lambda$.
\end{proof}

\begin{lem}
    \label{lem:lyapunov}
    The Lyapunov coefficient satisfies \(\Lbl = O(\lambda^{-1/6} \lvert n \rvert^{-1/3})\).
\end{lem}

\begin{proof}
Using Lemma~\ref{lem:covariance}, there exists a constant $C > 0$ such that
\[
\Lbl \leq C \sum_{x \in \XX} \left[\frac{\EEbl\lvert \overline{X_{1,x}}  \rvert^3 }{\lambda^{-1/2}}\frac{n_2^{1/2}}{n_1^{5/2}} + \frac{\EEbl\lvert \overline{X_{2,x}}  \rvert^3 }{\lambda^{-1/2}}\frac{n_1^{1/2}}{n_2^{5/2}} + \frac{\EEbl\lvert \overline{K_x}  \rvert^3 }{\lambda^{1/2}(n_1n_2)^{1/2}}\right].
\]
Therefore, we need only prove that
\[
 \sum_{x \in \L} \EEbl \left\lvert \overline{K_x}  \right\rvert^3 = O(\lvert n\rvert^{2/3}), \qquad
 \sum_{x \in \L} \EEbl \left\lvert \overline{X_{i,x}} \right\rvert^3 = O(\lvert n\rvert^{5/3}).
\]
Notice that for a Bernoulli random variable $B(p)$ of parameter $p$, one has $\EE[\lvert B(p) - p\rvert^3] \leq 4 (\EE [B(p)^3] + p^3) \leq 8 p$. This implies
\[
	 \sum_{x \in \L} \EEbl \left\lvert \overline{K_x}  \right\rvert^3 \leq \sum_{x \in \L} \frac{8 \lambda e^{-\beta\cdot x}}{1 - (1-\lambda) e^{-\beta\cdot x}} \leq \sum_{x \in \L} \frac{8 \lambda e^{-\beta\cdot x}}{1 - e^{-\beta\cdot x}} = O(\frac{\lambda}{\beta_1\beta_2}).
\]
Similarly, we obtain
\[	\sum_{x \in \L} \EEbl \left\lvert \overline{X_{1,x}} \right\rvert^3 = O(\frac{\lambda}{\beta_1^4\beta_2}),
	\quad
	\sum_{x \in \L} \EEbl \left\lvert \overline{X_{2,x}} \right\rvert^3 = O(\frac{\lambda}{\beta_1\beta_2^4}).
\]
\end{proof}

\begin{lem}
    \label{lem:cramer}
    Condition~\eqref{eq:cramer} of Lemma~\ref{lem:framework} is satisfied.
    More precisely,
    \[
        \limsup_{|n| \to +\infty} \quad\sup_{(t,u) \in [-\pi,\pi]^3\setminus\mathcal{E}_{\beta,\lambda}}\quad\frac{1}{\lambda^{1/3} |n|^{2/3}} \log |\phi_n(t,u)| < 0.
    \]
\end{lem}

\begin{proof}
    From Lemmas~\ref{lem:covariance} and~\ref{lem:lyapunov}, there exists a constant \(c > 0\) depending on \(\lambda\) such that for all \(n=(n_1,n_2)\) with \(|n|\) large enough,
    \[
        [-\pi,\pi]^3 \setminus \mathcal{E}_{\lambda,n} \subset
        \{(t,u) \in \RR^3 : c < |u| \leq \pi \text{ or } c \lambda^{1/3} |n|^{-1/3} < |t|\}.
    \]
    The strategy of the proof is to deal separately with the cases $|u| > c$ and $|t| > c \lambda^{1/3} |n|^{-1/3}$, which requires to find first adequate bounds for $|\phi_n(t,u)|$ in both cases. For all \((t_1,t_2, u) \in \RR^3\) and \(x\in \L\), let us write \(t = (t_1,t_2)\) and \(\rho^x = e^{-\beta \cdot x}\). The "partial" characteristic function \(\phi_n^x(t,u) = \EE[e^{i(t\cdot X_x + uK_x )}]\) is given by
    \[
        \phi_n^x(t,u) = \left(1 + \lambda e^{iu}\dfrac{e^{it\cdot x}\rho^x}{1-  e^{it\cdot x}\rho^x}\right)\left(1 + \lambda\dfrac{\rho^x}{1-  \rho^x}\right)^{-1},
    \]
    hence a straightforward calculation yields
    \begin{align*}
        \left|\phi_n^x(t,u)\right|^2 
                                     & = 1 - \frac{\frac{4\lambda\rho^x}{(1-(1-\lambda)\rho^x)^2}\left[\frac{\rho^x(2+(\lambda-2)\rho^x)}{(1-\rho^x)^2}|\sin(\frac{t\cdot x}{2})|^2 + |\sin(\frac{t\cdot x + u}{2})|^2 - \rho^x |\sin(\frac{u}{2})|^2 \right]}{1+\frac{4\rho^x}{(1-\rho^x)^2}|\sin(\frac{t\cdot x}{2})|^2}\\
                                     & \leq \exp\left\{ - \frac{\frac{4\lambda\rho^x}{(1-(1-\lambda)\rho^x)^2}\bigl(2\rho^x|\sin(\frac{t\cdot x}{2})|^2 + |\sin(\frac{t\cdot x + u}{2})|^2 - \rho^x |\sin(\frac{u}{2})|^2 \bigr)}{1+\frac{4\rho^x}{(1-\rho^x)^2}|\sin(\frac{t\cdot x}{2})|^2}\right\}
    \end{align*}
    Using the law of sines in a triangle with angles $\frac{t\cdot x}{2}$, $\frac{u}{2}$ and $\frac{2\pi - t\cdot x+u}{2}$, we see that the numerator inside the bracket is proportional (with positive constant) to
    \[
        2\rho^x \|a\|^2 + \|{b}\|^2 - \rho^x \|{a} + {b}\|^2
    \]
    where $a$ and $b$ are two-dimensional vectors. Since the real quadratic form
    $(a_i,b_i) \mapsto 2\rho\, a_i^2 + b_i^2 - \frac{2\rho}{1+2\rho} \,(a_i + b_i)^2$
    is positive for all $\rho \in (0,1)$ and for $i \in \{1,2\}$, we deduce that 
    \begin{equation}
        \label{eq:sinu}
        |\phi_n^x(t,u)| \leq
        \exp\left\{ - \frac{\frac{2\lambda\rho^x}{(1-(1-\lambda)\rho^x)^2}}{1+\frac{4\rho^x}{(1-\rho^x)^2}}\left(\frac{2\rho^x}{1+2\rho^x}-\rho^x\right)\left|\sin(\tfrac{u}{2})\right|^2 \right\}
    \end{equation}
    for all $x$ such that $\rho_x \leq \frac{1}{2}$.
    In the same way, the positivity of the quadratic form $(a_i,b_i)  \mapsto \frac{\rho}{1-\rho}\, a_i^2 + b_i^2 - \rho\, (a_i+b_i)^2$ yields
    \begin{equation}
        \label{eq:sintx}
        |\phi_n^x(t,u)| \leq
        \exp\left\{ - \frac{\frac{2\lambda\rho^x}{(1-(1-\lambda)\rho^x)^2}}{1+\frac{4\rho^x}{(1-\rho^x)^2}}\left(2\rho^x-\frac{\rho^x}{1-\rho^x}\right)\left|\sin(\tfrac{t\cdot x}{2})\right|^2 \right\}
    \end{equation}
    for all $x$ such that $\rho_x \leq \frac{1}{2}$.

    Let us begin with the region \(\{(t,u) \in \RR^3 : c < |u| \leq \pi\}\). In this case $|\sin(\tfrac{u}{2})|$ is uniformly bounded from below by $|\sin(\tfrac{c}{2})|$. Hence using \eqref{eq:sinu} for the $x \in \XX$ such that $\frac{1}{4} < \rho^x \leq \frac{1}{3}$ and the bound $|\phi_n^x(t,u)| \leq 1$ for all other $x$, we obtain
    \[
        \log |\phi_n(t,u)| \leq - \frac{1}{160}\frac{\lambda|\sin(\tfrac{c}{2})|^2}{(1+\frac{1}{3}|\lambda-1|)^2}\left|\left\{x \in \XX : \frac{1}{4} < \rho^x \leq \frac{1}{3}\right\}\right|.
    \]
    To conclude, let us recall that the number of integral points with coprime coordinates such that $\frac{1}{4} < e^{-\beta \cdot x} \leq \frac{1}{3}$ is asymptotically equal to $\frac{1}{\zeta(2)}\frac{\log(4/3)}{2\beta_1\beta_2} \asymp \lambda^{-2/3} |n|^{2/3}$.
    
    We now turn to the region \(\{(t,u) \in [-\pi,\pi]^3 : c\lambda^{1/3} |n|^{-1/3} < |t|\}\). Without loss of generality, we can assume \(|t_1| > c' \lambda^{1/3} |n|^{-1/3}\) for some universal constant \(c' \in (0; c)\). Using the inequality \eqref{eq:sintx} for the elements $x \in \XX$ such that $\frac{1}{4} < \rho^x \leq \frac{1}{3}$ and the bound $|\phi_n^x(t,u)| \leq 1$ for all other $x$, we obtain for all $\epsilon \in (0,1)$,
    \[
        \log |\phi_n(t,u)| \leq - \frac{\epsilon^2}{64} \frac{\lambda}{(1 + \frac{1}{3}|\lambda-1|)^2} \left|\left\{ x \in \XX : \frac{1}{4} < \rho^x \leq \frac{1}{3} \text{ and } |\sin(\tfrac{t \cdot x}{2})| \geq \epsilon\right\}\right|.
    \]
  Since the number of $x \in \XX$ such that $\frac{1}{4} < e^{-\beta \cdot x} \leq \frac{1}{3}$ is asymptotically equal to $\frac{\log(4/3)}{2\zeta(2)\beta_1\beta_2}$, it is enough to prove that we can find $\epsilon$ such that the set of vectors $x \in \ZZ_+^2$ with $|\sin(\tfrac{t\cdot x}{2})| < \epsilon$ has density strictly smaller than $\frac{1}{\zeta(2)}$ in $\{x \in \ZZ_+^2 : \frac{1}{4} < \rho^x \leq \frac{1}{3}\}$.
  We split up this region according to horizontal lines, that is to say with $\frac{t_2x_2}{2}$ constant. The set $\{x_1 \in \RR : |\sin(\frac{t_2x_2}{2}+\tfrac{t_1x_1}{2})| < \epsilon\}$ is a periodic union of strips of period $\tau_1 = \frac{2\pi}{t_1} \geq 2$ and width bounded by $4\epsilon\tau_1$. Hence the number of $x_1 \in \ZZ_+$ satisfying this condition and lying in any bounded finite interval $I$ is at most $\left(\frac{|I|}{\tau_1} + 2\right)(4\epsilon \tau_1+1)$. Summing up the contributions of the horizontal lines, this shows the existence of some positive constant $C > 0$ independent of $\epsilon$ such that for all $\epsilon \in (0,1)$, the number of $x \in \ZZ_+^2$ satisfying both $\frac{1}{4} < e^{-\beta \cdot x} \leq \frac{1}{3}$ and $|\sin(\tfrac{t\cdot x}{2})| < \epsilon$ is bounded by
  \[
      (\tfrac{1}{2} + C\epsilon) \frac{\log(4/3)}{2\beta_1\beta_2} + C |n|^{1/3}\log|n|.
  \]
  To achieve our goal, we can therefore choose $\epsilon = \frac{1}{2C}(\frac{1}{\zeta(2)}-\frac{1}{2}) > 0$.
\end{proof}

\begin{proof}[Proof of Theorem~\ref{thm:local_limit}]
    We simply check that the hypotheses of Lemma~\ref{lem:framework} are satisfied. From Lemma~\ref{lem:covariance}, we have $\sigmabl^2 \asymp k$ and $\det(\Gammabl) \asymp k^{-1}|n|^4$, hence
   \[
       \frac{1}{\sigmabl \sqrt{\det \Gammabl}} \asymp \frac{1}{|n|^2}.
   \]
   Using in addition Lemma~\ref{lem:lyapunov}, we have also
   \[
       \frac{\Lbl}{\sqrt{\det \Gammabl}} = O\left(\frac{1}{|n|^2}\right).
   \]
   Finally, Lemma~\ref{lem:cramer} shows the existence of some constant $c > 0$ such that for all $(n,k)$ large enough,
   \[
        \sup_{(t,u) \in [-\pi,\pi]^3
            \setminus\mathcal{E}_{\beta,\lambda}}
        \; |\phi_n(t,u)| \leq e^{-c k}
    \]
    Since we have made the assumption $\log |n| = o(k)$, the quantity $e^{-ck}$ is also bounded from above by $|n|^{-2}$. Therefore, all hypotheses of Lemma~\ref{lem:framework} are satisfied. As a consequence, $\PPbl$ satisfies a local limit theorem with speed rate $a_{\beta,\lambda} \asymp |n|^{-2}$. 

\end{proof}

\section{Limit shape}
\label{sec:limit_shape}

We start by proving the existence of a limit shape in the modified \Sinai{} model, which is the aim of the next two lemmas. The natural normalization for the convex chain is to divide each coordinate by the corresponding expectations for the final point.

The first lemma shows that the arc of parabola is the limiting curve of the expectation of the random convex chain $m_i^\theta (\beta,\lambda) = \EEbl [X_i^\theta]$ for $i\in \{1,2\}, \theta \in [0,\infty]$ under the $\PPbl$ distribution.

\begin{lem}
    \label{lem:expectation}
    Suppose that $\beta_1$ and $\beta_2$ tend to $0$ such that $\beta_1 \asymp \beta_2$ and $\lambda$ is bounded from above. Then
    \[
        \lim_{|\beta| \to 0} \sup_{\theta \in [0,\infty]} \left |
        \left[\frac{m_1^\theta(\beta,\lambda)}{m_1^\infty(\beta,\lambda)},
        \frac{m_2^\theta(\beta,\lambda)}{m_2^\infty(\beta,\lambda)}\right] - 
        \left[
            \frac{\theta(\theta + 2\frac{\beta_1}{\beta_2})}{(\theta + \frac{\beta_1}{\beta_2})^2} 
            ,
            \frac{\theta^2}{(\theta + \frac{\beta_1}{\beta_2})^2}
        \right] \right| = 0.
    \]
\end{lem}

\begin{proof}
    Since we are dealing with increasing functions, the uniform convergence convergence will follow from the simple convergence.
    We mimic the proof of Lemma~\ref{lem:derivatives_Z}, except that the domain of summation $\XX$ is replaced by the subset of vectors $x$ such that $x_2 \leq \theta x_1$. The expectations are given by the first derivatives of the \emph{partial} logarithmic partition function
    \[
        \log Z^\theta(\beta, \lambda) = \frac{1}{2i\pi}\int_{c-i\infty}^{c+i\infty} (\zeta(s+1) - \Li_{s+1}(1-\lambda)) {\zeta_2^{\theta,*}}(s)\Gamma(s)\,ds
    \]
    where $\zeta_2^{\theta,*}$ is the restricted zeta function defined by analytic continuation of the series
    \begin{align*}
        \zeta_2^{\theta,*}(s) & = \sum_{\substack{x \in \XX \\ x_2 \leq \theta x_1}} (\beta_1 x_1 + \beta_2 x_2)^{-s} \\
                              & = \frac{1}{\beta_1^s} + \frac{1_{\{\theta = \infty\}} }{\beta_2^s} + \frac{1}{\zeta(s)}\sum_{\substack{x_1,x_2 \geq 1\\ x_2 \leq \theta x_1}} (\beta_1 x_1 + \beta_2 x_2)^{-s}.
    \end{align*}

    The continuation of the underlying restricted Barnes zeta function is obtained using the Euler-Maclaurin formula several times:
    \begin{align*}
        \sum_{x_2 = 1}^{\lfloor \theta x_1 \rfloor} (\beta_1 x_1 + \beta_2 x_2)^{-s} & =
        \int_1^{\lfloor \theta x_1 \rfloor} (\beta_1 x_1 + \beta_2 x_2)^{-s}\,dx_2  + \frac{(\beta_1 x_1 + \beta_2)^{-s}}{2} + \frac{(\beta_1 x_1 + \beta_2 \lfloor \theta x_1 \rfloor)^{-s}}{2} \\
        & \qquad -s \beta_2 \int_1^{\lfloor \theta x_1 \rfloor}(\{x_2\} - \frac{1}{2}) (\beta_1 x_1 + \beta_2 x_2)^{-(s+1)}\,dx_2\\
        & = \int_1^{\theta x_1} (\beta_1 x_1 + \beta_2 x_2)^{-s}\,dx_2  + \frac{(\beta_1 x_1 + \beta_2)^{-s}}{2} + \frac{(\beta_1 x_1 + \beta_2 \lfloor \theta x_1 \rfloor)^{-s}}{2} \\
        & \qquad -s \beta_2 \int_1^{\lfloor \theta x_1 \rfloor}(\{x_2\} - \frac{1}{2}) (\beta_1 x_1 + \beta_2 x_2)^{-(s+1)}\,dx_2 \\
        & \qquad- \int_{\lfloor \theta x_1 \rfloor}^{\theta x_1} (\beta_1 x_1 + \beta_2 x_2)^{-s}\,dx_2\\
        & = \int_1^{\theta x_1} (\beta_1 x_1 + \beta_2 x_2)^{-s}\,dx_2  + \frac{(\beta_1 x_1 + \beta_2)^{-s}}{2} + \frac{(\beta_1 x_1 + \beta_2 \lfloor \theta x_1 \rfloor)^{-s}}{2} \\
        & \qquad -s \beta_2 \int_1^{\lfloor \theta x_1 \rfloor}(\{x_2\} - \frac{1}{2}) (\beta_1 x_1 + \beta_2 x_2)^{-(s+1)}\,dx_2 \\
        & \qquad- \int_{\lfloor \theta x_1 \rfloor}^{\theta x_1} (\beta_1 x_1 + \beta_2 x_2)^{-s}\,dx_2\\
        & = \frac{(\beta_1 x_1 + \beta_2)^{-s+1}}{\beta_2(s-1)} 
        - \frac{(\beta_1 x_1 + \beta_2 \theta x_1)^{-s+1}}{\beta_2(s-1)}  + R(s,x_1,\beta_1,\beta_2,\theta)
    \end{align*}
    where
    \begin{align*}
        R(s,x_1,\beta_1,\beta_2,\theta) & = \frac{(\beta_1 x_1 + \beta_2)^{-s}}{2} + \frac{(\beta_1 x_1 + \beta_2 \lfloor \theta x_1 \rfloor)^{-s}}{2} \\
        & \qquad -s \beta_2 \int_1^{\lfloor \theta x_1 \rfloor}(\{x_2\} - \frac{1}{2}) (\beta_1 x_1 + \beta_2 x_2)^{-(s+1)}\,dx_2 \\
        & \qquad- \int_{\lfloor \theta x_1 \rfloor}^{\theta x_1} (\beta_1 x_1 + \beta_2 x_2)^{-s}\,dx_2\\
    \end{align*}
    is such that $\sum_{x_1 \geq 1} R(s,x_1,\beta_1,\beta_2,\theta)$ converges absolutely for all $s$ with $\Re(s) > 1$. Therefore the latter series defines an holomorphic function in the half-plane $\Re(s) > 1$. Finally,

    \begin{align*}
        \sum_{\substack{x_1,x_2 \geq 1\\ x_2 \leq \theta x_1}} (\beta_1 x_1 + \beta_2 x_2)^{-s}
            & = \frac{(\beta_1+\beta_2)^{-s + 2}}{\beta_1\beta_2(s-1)(s-2)} - \frac{(\beta_1+\theta\beta_2)^{-s+2}}{(\beta_1+\theta\beta_2)\beta_2(s-1)(s-2)}\\
            & \qquad + \widetilde{R}(s,\beta_1,\beta_2,\theta)
    \end{align*}
    where $\widetilde{R}$ is holomorphic in $s$ for $\Re(s) > 1$. Hence, the residue at $s = 2$ is $\frac{\theta}{\beta_1(\beta_1 + \theta \beta_2)}$. 
Taking the derivatives with respect to $\beta_1$ and $\beta_2$, we obtain,
\begin{align*}
    - \frac{\partial}{\partial \beta_1}
    \sum_{\substack{x_1,x_2 \geq 1\\ x_2 \leq \theta x_1}} (\beta_1 x_1 + \beta_2 x_2)^{-s}
    & = \frac{1}{\beta_1^2\beta_2}\frac{\theta(\theta + 2 \frac{\beta_1}{\beta_2})}{(\theta + \frac{\beta_1}{\beta_2})^2} \frac{1}{s-2} + R_1(s,\beta_1,\beta_2,\theta)
\end{align*}
and similarly
\begin{align*}
    - \frac{\partial}{\partial \beta_2}
    \sum_{\substack{x_1,x_2 \geq 1\\ x_2 \leq \theta x_1}} (\beta_1 x_1 + \beta_2 x_2)^{-s}
    & = \frac{1}{\beta_1 \beta_2^2} \frac{\theta^2}{(\theta + \frac{\beta_1}{\beta_2})^2} \frac{1}{s-2} + R_2(s,\beta_1,\beta_2,\theta)
\end{align*}
where both remainder terms $R_1$ and $R_2$ are holomorphic in $s$ in the half-plane $\sigma := \Re(s) > 1$ and are bounded, up to positive constants, by
\[
    \frac{|s|^2}{\sigma - 1} \min(\beta_1,\beta_2)^{-\sigma-1}.
\]
This decrease makes it possible to apply the residue theorem in order to shift to the left the vertical line of integration from $\sigma = 3$ to $\sigma = \frac{3}{2}$.
When $\beta_1$ and $\beta_2$ tend to $0$ and $\frac{\beta_1}{\beta_2}$ tends to $\ell$, we thus find
\begin{align*}
    \EE_{\beta,\lambda}[X_1^\theta] &= \frac{\zeta(3)-\Li_3(1-\lambda)}{\zeta(2)}\left[ \frac{1}{\beta_1^2\beta_2} \frac{\theta(\theta + 2\frac{\beta_1}{\beta_2})}{(\theta + \frac{\beta_1}{\beta_2})^2} + O\left(\frac{1}{|\beta|^{5/2}}\right)\right],\\
    \EE_{\beta,\lambda}[X_2^\theta] &= \frac{\zeta(3)-\Li_3(1-\lambda)}{\zeta(2)} \left[\frac{1}{\beta_1^2\beta_2} \frac{\theta^2}{(\theta + \frac{\beta_1}{\beta_2})^2} + O\left(\frac{1}{|\beta|^{5/2}}\right)\right].
\end{align*}
We obtain the announced result by normalizing these quantities by their limits when $\theta$ goes to infinity.
\end{proof}

\begin{lem}[Uniform exponential concentration]
    \label{lem:concentration}
 Suppose that $\beta_1$ and $\beta_2$ tend to $0$ such that $\beta_1 \asymp \beta_2$ and $\lambda$ is bounded from above. For all $\eta \in (0,1)$, we have
    \[
        \PP_{\beta,\lambda}\left[\sup_{1\leq i \leq 2}\sup_{\theta \in [0,\infty]} \frac{|X^\theta_i - m_i^\theta(\beta,\lambda)|}{m_i^\infty(\beta,\lambda)} > \eta\right] \leq
        \exp\left\{-\frac{c(\lambda)\eta^2}{8\beta_1\beta_2}\left(1 + o(1)\right)\right\}.
    \]
\end{lem}

\begin{proof}
    Fix $i \in \{1,2\}$ and let $M_\theta = X_i^\theta - m_i^\theta(\beta,\lambda)$ for all $\theta \geq 0$. The stochastic process $(M_\theta)_{\theta \geq 0}$ is a $\PPbl$-martingale, therefore $(e^{t M_\theta})_{\theta \geq 0}$ is a positive $\PPbl$-submartingale for any choice of $t \geq 0$ such that $\EEbl[e^{tX_i}]$ is finite. This condition is satisfied when $t < \beta_1$. Doob's martingale inequality implies for all $\eta > 0$,
    \begin{align*}
        \PPbl\left[\sup_{\theta \in [0,\infty]} M_\theta > \eta\,m_i^\infty(\beta,\lambda)\right]
        & = \PPbl\left[\sup_{\theta \in [0,\infty]} e^{tM_\theta}
            > e^{t\eta m_i^\infty(\beta,\lambda)}\right]\\
        & \leq e^{-t\eta m_i^\infty(\beta,\lambda)} \,\EEbl\left[e^{tM_\infty}\right] = e^{-t(\eta+1)m_i^\infty(\beta,\lambda)} \,\EEbl[e^{tX_i}]
    \end{align*}
    For $i = 1$, Lemma~\ref{lem:derivatives_Z} shows that the logarithm of the right-hand side satisfies
    \[
        -t(1+\eta) m_1^\infty(\beta,\lambda) + \log \frac{Z(\beta_1 - t, \beta_2 ; \lambda)}{Z (\beta_1,\beta_2;\lambda)} = \frac{c(\lambda)}{\beta_1\beta_2}\left[-\frac{t(1+\eta)}{\beta_1} - 1 + \frac{\beta_1}{\beta_1-t} + o(1)\right]
    \]
    asymptotically when $t$ and $\beta_1$ are of the same order. The same holds for $i=2$. This is roughly optimized for the choice $t = \beta_i\left(1-(1+\eta)^{-1/2}\right)$, which gives
    \[
        \PP_{\beta,\lambda} \left[\sup_{\theta \in[0,\infty]} M_\theta > \eta\, m_i^\infty(\beta,\lambda)\right] \leq \exp\left\{-\frac{2c(\lambda)}{\beta_1\beta_2}\left(1+\frac{\eta}{2}-\sqrt{1+\eta} + o(1)\right)\right\}.
    \]
    When considering the martingale defined by $N_\theta = m_i^\theta(\beta,\lambda) - X_i^\theta$, one obtains with the same method
    \[
        \PPbl\left[\sup_{\theta \in [0,\infty]} N_\theta > \eta\,m_i^\infty(\beta,\lambda)\right] \leq \exp\left\{-\frac{2c(\lambda)}{\beta_1\beta_2}\left(1 - \frac{\eta}{2} - \sqrt{1-\eta} + o(1)\right)\right\}.
    \]
    Since the previous inequalities hold for both $i\in\{1,2\}$, a simple union bound now yields
    \[
        \PPbl\left[\sup_{1\leq i \leq 2}\sup_{\theta \in [0,\infty]} \frac{|X^\theta_i - m_i^\theta(\beta,\lambda)|}{m_i^\infty(\beta,\lambda)} > \eta\right] \leq 4 \exp\left\{-\frac{c(\lambda)\eta^2}{8\beta_1\beta_2}\left(1+ o(1)\right)\right\}.
    \]
\end{proof}

We introduce the following parametrization of the arc of parabola $\sqrt{y} + \sqrt{1-x} = 1$:
\[
    x_1(\theta) = \frac{\theta(\theta+2)}{(\theta + 1)^2}, \quad x_2(\theta) = \frac{\theta^2}{(\theta + 1)^2}, \qquad \theta \in [0,\infty].
\]

\begin{theorem}[Limit shape for numerous vertices]
	\label{thm:limit_shape_numerous}
    Assume that $n_1 \asymp n_2 \to +\infty$, and $k = O(|n|^{2/3})$, and $\log |n| = o(k)$.
   There exists $c > 0$ such that for all $\eta \in (0,1)$, 
    \[
        \PP_{n,k}\left[\sup_{1\leq i \leq 2}\sup_{\theta \in [0,\infty]} \frac{|X^\theta_i - x_i(\frac{\beta_2}{\beta_1}\theta)|}{n_i} > \eta\right] \leq
        \exp\left\{-c\eta^2 k\left(1+ o(1)\right)\right\}.
    \]
    In particular,  the Hausdorff distance between a random convex chain on $\frac{1}{n}\ZZ_+^2$ joining $(0,0)$ to $(1,1)$ with at most $k$ vertices and the arc of parabola $\sqrt{y}+ \sqrt{1-x} = 1$ converges in probability to $0$.
\end{theorem}

\begin{proof}
    Using the triangle inequality and Lemma~\ref{lem:expectation}, we need only prove the analogue of Lemma~\ref{lem:concentration} for the uniform probability $\PP_{n,k}$. Remind that the measure $\PPbl$ conditional on the event $\{X=n,K=k\}$ is nothing but the uniform probability $\PP_{n,k}$. Hence for all event $E$,
    \[
        \PP_{n,k}(E) \leq \frac{\PPbl(E)}{\PPbl(X=n,K=k)}.
    \]
    Applying this with the deviation event above for the parameters $(\beta,\lambda)$ defined in Section~\ref{sec:calibration} and using the Local Limit Theorem~\ref{thm:local_limit} as well as the concentration bound provided by Lemma~\ref{lem:concentration}, the right-hand side reads, up to constants,
    \[
        \frac{|n|^2}{\sqrt{k}} \,\exp\left\{-c\eta^2 k (1+o(1))\right\}.
    \]
    Since $\log |n| = o(k)$, the result follows.
\end{proof}

\input{very_few_vertices}

\section{Back to \Jarnik's problem}
\label{sec:jarnik}

In \cite{MR1544776}, \Jarnik gives an asymptotic formula of the maximum possible number of vertices of a convex lattice polygonal line having a \emph{total Euclidean length} smaller than $n$, and whose segments make an angle with the $x$-axis between $0$ and $\frac{\pi}{4}$. What he finds is $\frac{3}{2}\,\frac{n^{2/3}}{(2\pi)^{1/3}}$. If, in order to be closer to our setting, we  ask the segments to make an angle with the $x$-axis between $0$ and $\frac{\pi}{2}$, \Jarnik's formula is
changed into $\frac{3}{2}\frac{n^{2/3}}{\pi^{1/3}}$ (which is  twice the above result for $\frac{n}{2}$). 

In this section, we want to present a detailed combinatorial analysis of this set of lines, which leads to  \Jarnik's result as well as to the asymptotic of the \emph{typical} number of vertices of such lines. It is the analog
of \Barany, \Sinai and Vershik's result when the constraint concerns the total length.

Let us first describe \Jarnik's argument, which is a good application of the correspondence described in Section~\ref{sec:correspondence}. It says the following: the function $\omega$ realizing the maximum can be taken among the functions taking their values in $\{0,1\}$. Indeed, by changing the non-zero values of a function $\nu$ into $1$, one can obtain a chain with the same number of vertices, but with a shorter length. Now, if the number of vertices $k$ is given, the convex chain having minimal
length, will be defined by the function $\omega$ which associates $1$ to the $k$ points  of $\XX$ which are the closest to the origin.
Since the set $X$ has an asymptotic density $\frac{6}{\pi^2}$, when  $N$ is big, this set  of points is asymptotically equivalent to the intersection of $X$ with the disc of center $O$ having radius $R$ satisfying $\frac{6}{\pi^2}\cdot  \frac{\pi R^2}{4}=N$ i.e. $R=(\frac{2\pi}{3} N)^{1/2}$. The total length of the line is equivalent to 
$L= \int_0^R r\times \frac{6}{\pi^2}\frac{\pi}{2}r dr=\frac{R^3}{\pi}=\frac{(\frac{2\pi}{3} N)^{3/2}}{\pi}$. 
This yields precisely $N=\frac{3}{2} \frac{L^{2/3}}{\pi^{1/3}}\simeq1.02\, L^{2/3}$. 

In order to get finer results, we introduce the probability distribution on the space $\Omega$ proportional to
\[
    \exp\left(-\beta \sum_{x \in \XX} \omega(x) \sqrt{|x_1|^2 + |x_2|^2}\right) \lambda^{\sum_{x\in \XX} 1_{\{\omega(x) > 0\}}}
\]
which depends on two parameters $\beta,\lambda$. In this set-up, the partition function turns out to be
\[
    Z = \prod_{x\in \XX} \frac{1-(1-\lambda)e^{-\beta \sqrt{|x_1|^2 + |x_2|^2}}}{1-e^{-\beta \sqrt{|x_1|^2+|x_2|^2}}}.
\]
The Mellin transform representation for $\log Z$ now involves
\[
    \frac{\Gamma(s)(\Li_{s+1}(1-\lambda)-\zeta(s+1))}{\zeta(s)}\sum_{x_1,x_2 \geq 1} {(|x_1|^2+|x_2|^2)}^{-s/2}, \qquad \Re(s) > 2.
\]
The factors $\zeta(s)^{-1}$ and $\Li_{s+1}(1-\lambda)-\zeta(s+1)$, which correspond respectively to the coprimality condition on the lattice and to the penalization of vertices, are still present. The main difference relies in the replacement of the Barnes zeta function by the Epstein zeta function which comes from the penalization by length in the model. With the help of the residue analysis of this Mellin transform and a local limit theorem, we obtain:

\begin{theorem}
    Let $p_J(n;k)$ denote the number of convex chains on $\ZZ_+^2$ issuing from $(0,0)$ with $k$ vertices and length between $n$ and $n+1$. As $n$ tends to $+\infty$,
    \[
        \text{if}\qquad \frac{k}{n^{2/3}} \longrightarrow \frac{\pi^{1/3}}{2}\c(\lambda), \qquad\text{then}\qquad  \frac{1}{n^{2/3}} \log p_J(n;k) \longrightarrow \frac{\pi^{1/3}}{2}\e(\lambda),
    \]
    where $\e$ and $\c$ are the functions introduced in Theorem~\ref{thm:detailed_comb}. Moreover, the Hausdorff distance between a random element of this set normalized by $\frac{1}{n}$, and the arc of circle $\{(x,y) \in [0,1]^2 : x^2 + (y-1)^2 = 1\}$ converges to $0$ in probability.
\end{theorem}

From this result, we deduce that the typical number of vertices of such a chain which is achieved for $\lambda=1$ is asymptotically equal to
\[
    \left(\frac{3}{4\pi\zeta(3)^2}\right)^{1/3} n^{2/3}.
\]
Similarly, the total number of convex chains having length between $n$ and $n+1$ is asymptotically equal to
\[
    \exp\left(\frac{3^{4/3}\zeta(3)^{1/3}}{(4\pi)^{1/3}}\, n^{2/3} (1+o(1))\right).
\]
In addition, we can derive \Jarnik's result in the lines of Remark~\ref{rem:longest}.

\section{Mixing constraints and finding new limit shapes}
\label{sec:onion}

In this section we introduce a family of lattice convex chain models which achieves a continuous interpolation of limit shapes between the diagonal of the square and the South-East corner sides of the square, passing through the arc of circle and the arc of parabola. Let $\|\cdot\|_1$ and $\|\cdot\|_2$ denote respectively the Taxicab norm and the Euclidean norm on $\RR^2$. Recall that for all $x \in \RR^2$,
\[
    \|x\|_1 = |x_1| + |x_2| \geq \|x\|_2 = \sqrt{|x_1|^2+|x_2|^2} \geq \frac{1}{\sqrt{2}}\|x\|_1.
\]

The Gibbs distribution we consider on the space $\Omega$ involves both these norms in order to take into account both the extreme point of the chain and its length:
\[
    \frac{1}{Z}\exp\left(-\beta \sum_{x\in\XX} \omega(x) (\|x\|_1 + \lambda\sqrt{2} \|x\|_2)\right), \quad
    Z = \prod_{x \in \XX} \left(1 - e^{-\beta(\|x\|_1 + \lambda\sqrt{2}\|x\|_2)}\right).
\]
This infinite product is convergent if $\beta >0$ and $\lambda > -\frac{1}{\sqrt{2}}$ or if $\beta < 0$ and $\lambda < -1$.
In both cases, the Mellin transform representation of $\log Z$ involves
\[
    \frac{\Gamma(s)\zeta(s+1)}{\zeta(s)} \sum_{x_1,x_2 \geq 1} (\|x\|_1 + \lambda\sqrt{2} \|x\|_2)^{-s},\qquad \Re(s) > 2.
\]
As usual, the leading term of the expansion of $\log Z$ when $\beta \to 0$ is obtained by computing the residue of this function at $s = 2$. It turns out to be
\[
    \frac{\zeta(3)}{2\zeta(2)} \int_{-\pi/4}^{\pi/4} \frac{d\theta}{(\lambda + \cos (\theta))^2}.
\]
An application of the residue theorem shows that the expected length of the curve is asymptotically equivalent to
\[
    \frac{1}{\beta^3}\frac{\zeta(3)}{\sqrt{2}\zeta(2)} \int_{-\pi/4}^{\pi/4} \frac{d\theta}{(\lambda + \cos (\theta))^3}
\]
and that the coordinates of the ending point have asymptotic expected value
\[
    \frac{1}{\beta^3}\frac{\zeta(3)}{2\zeta(2)} \int_{-\pi/4}^{\pi/4} \frac{\cos(\theta) d\theta}{(\lambda + \cos (\theta))^3}.
\]

\begin{figure}[h]
\begin{center}
    \includegraphics[scale=1.2]{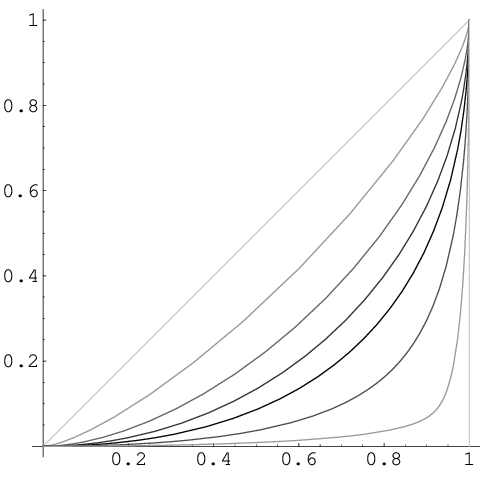}
    \caption{Limit shapes of different Euclidean lengths. Successively: $\sqrt{2}$ (diagonal); $1.454$, $1.516$, $\frac{\pi}{2}$ (circle), $1+\frac{\ln(1+\sqrt{2})}{\sqrt{2}}$ (parabola), $1.716$, $1.861$ and $2$ (square).}
    \label{fig:onion}
\end{center}
\end{figure}

As in previous sections, a local limit theorem gives a correspondence between this Gibbs measure and the uniform distribution on a specific set of convex chains, namely the convex chains with endpoint $(n,n)$ and total length belonging to $[L\cdot n, L\cdot n + 1]$ for some $L \in ]\sqrt2, 2[$ which is a function of $\lambda$,
\[
    L(\lambda) = \sqrt2{\int_0^{\pi\over4}{1\over(\lambda+\cos u)^3}du\over\int_0^{\pi\over4}{\cos u\over(\lambda+\cos u)^3}du}.
\]
By computations analogous to Section~\ref{sec:limit_shape}, one can show that the uniform distribution of chains with length between $L(\lambda)\cdot n$ and
$L(\lambda) \cdot n+1$ concentrates around the curve described by the parametrization
$$x_\lambda(\phi)=\sqrt2{\int_{0}^\phi{\cos u\over(\lambda+\cos (u-{\pi\over4}))^3}du\over\int_{-\pi/4}^{\pi/4}{\cos u\over(\lambda+\cos u)^3}du},\quad y_\lambda(\phi)=\sqrt2{\int_{0}^\phi{\sin u\over(\lambda+\cos (u-{\pi\over4}))^3}du\over\int_{-\pi/4}^{\pi/4}{\cos u\over(\lambda+\cos u)^3}du} \quad (0\leq\phi\leq{\pi\over2}).$$
The table provided in Figure~\ref{fig:table} resumes the limit shapes that we obtain for some limit values of $\lambda$. See also Figure~\ref{fig:onion} for a plot showing the interpolation of those limit shapes.

\begin{figure}[h]
\begin{center}
    \renewcommand{\arraystretch}{2}
    \begin{tabular}{|c|cc||ccc|}
        \hline
        $\lambda$ & $-\infty$ & $-1$ & $-\dfrac{1}{\sqrt{2}}$ & $0$ & $+\infty$\\
        \hline
        Limit shape & circle & diagonal & square & parabola & circle\\
        Length $L(\lambda)$ & $\dfrac{\pi}{2}$ & $\sqrt{2}$ & $2$ & $1+\dfrac{\ln(1+\sqrt{2})}{\sqrt{2}}$ & $\dfrac{\pi}{2}$\\
        \hline
    \end{tabular}
\end{center}
\caption{Critical and special values in the spectrum of limit shapes for the model of lattice convex chains with mixed constraints.}
\label{fig:table}
\end{figure}

\section{Towards an asymptotic equivalent for the number of convex chains}
\label{sec:sf}

In this section, we \emph{formally} push some steps further the asymptotic analysis of Section~\ref{sec:dca} under the assumption of a strong form of the Riemann Hypothesis, namely that all non trivial zeros of $\zeta(s)$ lie on the vertical line $\Re(s) = \frac{1}{2}$ and are simple.

In order to make the calculations easier to follow, we specify the model of Section~\ref{sec:dca} with $\lambda = 1$ and $\beta_1 = \beta_2 = \beta$, which corresponds to \Sinai's original model. In this case, the formula of Lemma~\ref{lem:integral_partition_function} for the logarithmic partition function may be written
\[
    \log Z = \frac{1}{2i\pi} \int_{c-i\infty}^{c+i\infty} \frac{\Gamma(s) \zeta(s+1)(\zeta(s-1)+\zeta(s))}{\zeta(s)\beta^s}\,ds
\]
for all $\beta > 0$ and $c > 2$. As in Lemma~\ref{lem:derivatives_Z}, the residue theorem yields the following formal asymptotic expansion where the sum is taken over all zeros $\rho$ of $\zeta(s)$ with $\Re(\rho) = \frac{1}{2}$,
\begin{equation}
    \label{eq:sf_log}
    \log Z = \frac{\zeta(3)}{\zeta(2)}\frac{1}{\beta^2} + \sum_{\rho} \frac{\Gamma(\rho)\zeta(\rho+1)\zeta(\rho-1)}{\zeta'(\rho) \beta^\rho} + \frac{7}{6} \log \frac{1}{\beta} + C  + o(1), \qquad \beta \to 0
\end{equation}
with $C = -2 \zeta'(-1) - \frac{1}{6} \log(2\pi)$. By considering the derivative of this expansion, we see that the calibrated parameter $\beta$ which is defined by $\frac{d}{d\beta} \log Z = - 2n$ satisfies
\begin{equation}
    \label{eq:sf_beta}
    \frac{1}{\beta^3} = \frac{n}{\kappa} - \sum_\rho \frac{\Gamma(\rho+1)\zeta(\rho+1)\zeta(\rho-1)}{\zeta'(\rho)} \left(\frac{n}{\kappa}\right)^{\frac{\rho+1}{3}} - \frac{7}{12 \kappa} \left(\frac{n}{\kappa}\right)^{1/3} + O(1).
\end{equation}
with $\kappa = \zeta(3)/\zeta(2)$. For this calibrated parameter, a local limit theorem similar to Theorem~\ref{thm:local_limit} gives
\begin{equation}
    \label{eq:sf_llt}
    \frac{1}{Z} e^{-2n \beta} p(n) = \PP_\beta[X=n] \sim \frac{\kappa^{1/3}}{2\pi \sqrt{3} n^{4/3}}
\end{equation}
where the right-hand term corresponds, up to $2\pi$, to the inverse of the standard deviation of $X$ under $\PP_\beta$. This is nothing but the square root of the determinant of the Hessian matrix of $\log Z$ at $\beta$ which is asymptotically equal to
\[
\begin{vmatrix}\frac{2\kappa}{\beta^4} & \frac{\kappa}{\beta^4}\\ \frac{\kappa}{\beta^4} & \frac{2\kappa}{\beta^4}\end{vmatrix} = \frac{3\kappa^2}{\beta^8}.
\]
Gathering \eqref{eq:sf_log}, \eqref{eq:sf_beta}, \eqref{eq:sf_llt}, we obtain the following asymptotic equivalent for the number $p(n)$ of lattice convex chain in $[0,n]^2$ from $(0,0)$ to $(n,n)$:
\[
    p(n) \sim \frac{e^{-2\zeta'(-1)}}{(2\pi)^{7/6} \sqrt{3} \kappa^{1/18} n^{17/18}}\exp\left[3\kappa^{1/3} n^{2/3} + \sum_\rho \frac{\Gamma(\rho)\zeta(\rho+1)\zeta(\rho-1)}{\zeta'(\rho)} \left(\frac{n}{\kappa}\right)^{\rho/3}\right].
\]

\bibliography{references}

\end{document}

%% file: very_few_vertices.tex
\section{Chains with few vertices}
\label{sec:few}

\subsection{Combinatorial analysis}
The previous machinery does not apply in the case of very few vertices but it can be completed by an an elementary approach that we present now which will actually work up to a number of vertices negligible compared to $n^{1/3}$. It is based on the following heuristics: when $n$ tends to $+\infty$ and the number of edges $k$ is very small compared to $n$, one can expect that choosing an element of $\Pi(n;k)$ at random is somewhat similar to choosing $k-1$ vertices from $[0,1]^2$ in convex position at random. 
\Barany~\cite{barany_sylvester_1999} and \Barany, Rote, Steiger, Zhang~\cite{barany_al_central_2000} proved by two different methods the existence of a parabolic limit shape in this continuous setting. These works are based on Valtr's observation that each convex chain with $k$ edges is associated, by permutation of the edges, to exactly $k!$ increasing North-East polygonal chains with pairwise different slopes.
%
%
%

Our first theorem is the convex-chain analogue to a result of Erd{\"o}s and Lehner on integer partitions \cite{erdos_lehner_distribution_1941}*{Theorem~4.1}.
\begin{theorem}
\label{thm:asymp_very_few_vertices}
    The number of convex chains joining $(0,0)$ to $(n,n)$ with $k$ edges satisfies
    \[
        \left|\Pi(n;k)\right| = \frac{1}{k!}\binom{n-1}{k-1}^2\left(1+o(1))\right),
    \]
    this formula being valid uniformly in $k$ for $k = o(n^{1/2}/(\log n)^{1/4})$.

\end{theorem}

\begin{proof}
    Let us start by proving an upper bound. This is done by considering the inequality
    \[
        \left| \Pi(n;k) \right| \leq \frac{1}{k!}\binom{n-1}{k-1}^2 + \frac{2}{(k-1)!}\binom{n-1}{k-2}\binom{n-1}{k-1} + \frac{1}{(k-2)!} \binom{n-1}{k-2}^2
    \]
    where the first term bounds the number of convex chains which are associated to strictly North-East chains, the second term bounds the number of convex chains having either a first horizontal vector or a last vertical one, and the third term bounds the numbers of convex chains having both a horizontal and a vertical vector. 

We now turn to a lower bound. Let $\{U_1, U_2, \dots, U_{k-1}\}$ and $\{V_1, V_2, \dots, V_{k-1}\}$ be two independent uniformly random subsets of $\{1,\dots,n-1\}$ of size $k-1$ whose elements are indexed in increasing order $U_1 < U_2 < \cdots < U_{k-1}$ and $V_1 < V_2 < \cdots <  V_{k-1}$. Let $M_0 = (0,0)$, $M_k = (n,n)$ and $M_i = (U_i,V_i)$ for $1 \leq i \leq k-1$. Obviously, the polygonal chain $(M_0,M_1,\dots,M_n)$ has uniform distribution among all increasing polygonal chain from $(0,0)$ to $(n,n)$.
We claim that the distribution of
$(\overrightarrow{M_{0}M_1},\overrightarrow{M_1M_2},\dots,\overrightarrow{M_{k-1}M_k})$ conditioned on the event that no two of these vectors are parallel is uniform among the chains of $\Pi(n,k)$ such that no side is parallel to the $x$-axis or the $y$-axis. Moreover, since the vectors are exchangeable, the probability that we can find $i < j$ such that $\overrightarrow{M_{i-1}M_i}$ and $\overrightarrow{M_{j-1}M_j}$ are parallel is bounded from above by $\binom{k}{2}$ times the probability that
$Y = \overrightarrow{M_0M_1}$ and $Z = \overrightarrow{M_1M_2}$ are parallel. Using the simple estimate
\[
    \binom{n-1}{k-1} \geq \frac{n^{k-1}}{(k-1)!}(1-o(1))
\]
which is asymptotically true since $k = o(\sqrt{n})$, we find that for all $(y,z) \in (\N^2)^2$, the probability that $Y=y$ and $Z=y$ is
\begin{align*}
    \PP(Y = y, Z=z) &= \frac{\binom{n-y_1-z_1}{k-3}\binom{n-y_2-z_2}{k-3}}{\binom{n-1}{k-1}^2}\\
    & \leq \frac{4k^2}{n^2}\left(1 - \frac{y_1+z_1}{n}\right)_+^{k-3}\left(1-\frac{y_2+z_2}{n}\right)_+^{k-3}\\
    & \leq \frac{4k^2}{n^2} \exp\left\{-\frac{k-3}{n}\left(y_1+y_2+z_1+z_2\right)\right\}.
\end{align*}
We can therefore dominate the probability that $Y$ and $Z$ are parallel by the probability that geometrically distributed random vectors are parallel, which is exactly estimated in the following lemma applied with $\beta = \frac{k}{n}$. In conclusion, the probability that at least two vectors are parallel is bounded by $\frac{k^4}{n^2}\log(n)$ up to a constant.
\end{proof}

\begin{lem}
    Let $Y_1,Y_2,Z_1,Z_2$ be independent and identically distributed geometric random variables of parameter $1-e^{\beta}$ with $\beta > 0$. When $\beta$ goes to $0$, the probability that the vectors $Y=(Y_1,Y_2)$ and $Z = (Z_1,Z_2)$ are parallel is asymptotically equal to
    \[
        \frac{\beta^2}{\zeta(2)}\log \frac{1}{\beta}.
    \]
\end{lem}

\begin{proof}
    The probability that $Y$ and $Z$ are parallel is
    \[
        \sum_{x\in \XX}\sum_{i,j \geq 1} \PP(Y = i\,x, Z = j \,x) = (1-e^{-\beta})^4 \sum_{x\in\XX}\sum_{i,j \geq 1} e^{-\beta(i+j)(x_1 + x_2)}.
    \]
    The Mellin transform of the double summation in the right-hand side with respect to $\beta > 0$ is well-defined for all $s \in \CC$ with $\Re(s) > 2$ and it is equal to
    \[
        \sum_{x\in \XX}\sum_{i,j \geq 1} \frac{\Gamma(s)}{(x_1+x_2)^s(i+j)^s} = \frac{\Gamma(s)}{\zeta(s)}(\zeta(s-1)-\zeta(s))^2.
    \]
    Expanding this Mellin transform in Laurent series at the pole $s = 2$ of order $2$ and using the residue theorem to express the Mellin inverse, one finds
    \[
        \sum_{x\in \XX}\sum_{i,j \geq 1} e^{-\beta(i+j)(x_1+x_2)} = \frac{1}{\zeta(2)}\frac{\log \frac{1}{\beta}}{\beta^2} - \frac{C}{\beta^2} + O\left(\frac{1}{\beta}\right),\qquad \text{as }\beta\to 0.
    \]
    where $C = \frac{2\zeta(2) - \zeta'(2) - 1 - \gamma}{\zeta(2)} \approx 0.471207$.
\end{proof}

\subsection{Limit shape}

\begin{theorem}[Limit shape for few vertices]
\label{thm:limit_shape_few}
The Hausdorff distance between a random convex chain in $(\frac{1}{n}\ZZ \cap [0,1])^2$ joining $(0,0)$ to $(1,1)$ having at most $k$ vertices and the arc of parabola $\sqrt{y}+ \sqrt{1-x} = 1$ converges in probability to $0$ when both $n$ and $k$ tend to $+\infty$ with $k = o(n^{1/3})$.
\end{theorem}

\begin{proof}
\Barany~\cite{barany_sylvester_1999} and \Barany, Rote, Steiger, Zhang~\cite{barany_al_central_2000} proved by two different methods the existence of a limit shape in the following continuous setting: if one picks at random $k-1$ points uniformly from the square $[0,1]^2$, then conditional on the event that these points are in convex position, the Hausdorff distance between the convex polygonal chain thus defined and the parabolic arc goes to $0$ in probability as $k$ goes to $+\infty$. Our strategy is to show that this result can be extended to the discrete
setting $([0,1] \cap \frac{1}{n}\ZZ)^2$ if $k$ is small enough compared to $n$ by using a natural embedding of the discrete model into the continuous model.

For this purpose, we first observe that the distribution of the above continuous model can be described as follows: pick uniformly at random $k-1$ points from both the $x$-axis and the $y$-axis, rank them in increasing order and let $0 = U_0  < U_1 <
U_2 < \dots < U_{k-1} < U_k = 1$ and $0 = V_0 < V_1 < V_2 < \dots < V_{k-1} < V_k = 1$ denote this ranking. The points $(U_i,V_i)$ define an increasing North-East polygonal chain joining $(0,0)$ to $(1,1)$. Reordering the segment lines of this chain by increasing slope order, exchangeability arguments show that we obtain a convex chain with $k$ edges that follows the desired distribution. This is analogous to the discrete construction of strictly North-East convex chains from $(0,0)$ to $(n,n)$ that occurs in the proof of Theorem~\ref{thm:asymp_very_few_vertices}.

Now, we define the lattice-valued random variables $\tilde{U}_0 \leq \tilde{U}_1 \leq \tilde{U}_2 \leq \dots \leq \tilde{U}_{k-1} \leq \tilde{U}_{k}$ and  $\tilde{V}_0  \leq \tilde{V}_1 \leq \tilde{V}_2 \leq \dots \leq \tilde{V}_{k-1} \leq \tilde{V}_k$ by discrete approximation:
\[
    \begin{cases}
        \tilde{U}_i \in \frac{1}{n}\ZZ, \quad U_i \leq \tilde{U}_i < U_i + \frac{1}{n}\\
        \tilde{V}_i \in \frac{1}{n}\ZZ, \quad V_i - \frac{1}{n} < \tilde{V}_i \leq V_i,
    \end{cases}\qquad \text{for }1 \leq i \leq k-1.
\]
Remark that we still have $(\tilde{U}_0,\tilde{V}_0) = (0,0)$ and $(\tilde{U}_k,\tilde{V}_k) = (1,1)$.

Let $X_i = (U_i - U_{i-1} ,V_i - V_{i-1})$ and let $\tilde{X}_i = (\tilde{U}_i - \tilde{U}_{i-1} ,\tilde{V}_i - \tilde{V}_{i-1})$ be the discrete approximation of $X_i$ for $1 \leq i \leq k$. Conditional on the event that the slopes of $(X_1,\dots,X_k)$ and $(\tilde{X}_1,\dots, \tilde{X}_k)$ are pairwise distinct and ranked in the same order, the Hausdorff distance between the associated convex chains is bounded by $\frac{k}{n}$, which goes asymptotically to $0$. Since a direct
application of \cite{barany_al_central_2000}*{Theorem 2} shows that the distance between the convex chain associated to $X$ and the parabolic arc converges to $0$ in probability as $k$ tends to $+\infty$, we deduce
that the Hausdorff distance between the convex chain associated to $\tilde{X}$ and the parabolic arc also converges in probability to $0$ on this event.
As in the proof of Theorem~\ref{thm:asymp_very_few_vertices}, the joint density of $(X_i,X_j)$ is dominated by the density of a couple of independent vectors
whose coordinates are independent exponential variables with parameter $k$. These vectors being of order of magnitude $\frac{1}{k}$, the order of the slopes of $(X_i,X_j)$ and $(\tilde{X}_i,\tilde{X}_j)$ may be reversed only if the angle between $X_i$ and $X_j$ is smaller than $\frac{ck}{n}$ for some $c > 0$, which happens with probability of order $\frac{k}{n}$. Henceforth, the probability that there exists $i < j$ for which the slopes of $(X_i,X_j)$ and $(\tilde{X}_i,\tilde{X}_j)$ are ranked in opposite is bounded, up to a constant, by $\binom{k}{2} \frac{k}{n}$. Therefore, the Hausdorff distance between the convex chain associated to $\tilde{X}$ and the parabolic arc also converges to $0$ in probability if $k = o(n^{1/3})$.

The final step is to compare the distribution of the increasing reordering of $(\tilde{X}_1,\dots,\tilde{X}_k)$ with the uniform distribution on $\Pi(n;k)$. As a consequence of Theorem~\ref{thm:asymp_very_few_vertices}, the probability that a uniformly random element of $\Pi(n;k)$ is strictly North-East tends to $1$. The key points, which follows from Valtr's observation, is that the uniform distribution on strictly North-East convex chains with $k$ edges coincides with the distribution of the chain obtained by reordering the vectors $(\tilde{X}_1,\dots,\tilde{X}_k)$, conditional on the event that these vectors are pairwise linearly independent and strictly North-East. Since we showed in the previous paragraph that all the angles between two vectors of $(\tilde{X}_1,\dots,\tilde{X}_k)$ are at least $\frac{ck}{n}$ with probability $1 - O(\frac{k^3}{n})$, the linear independence condition occurs with probability tending to $1$.  On the other hand, $(\tilde{X}_1,\dots,\tilde{X}_k)$ are strictly North-East with probability $1 - O(\frac{k^2}{n})$. Therefore, the event we conditioned on has a probability tending to $1$, which proves that the total variation distance between the two distributions tends to $0$.
\end{proof}